\documentclass[smallextended,envcountsect]{svjour3} 
\smartqed 
\usepackage{graphicx}
\journalname{}
\usepackage{amsfonts}
\usepackage{rotating}
\usepackage{lscape}
\usepackage{supertabular}
\usepackage{longtable}
\usepackage{pst-plot}
\usepackage{pstricks}
\usepackage{pst-text}
\begin{document}
\title{
\bf An affine scaling method using a class of differential barrier functions}\date{ }
\author{Abdessamad BARBARA}
\institute{Abdessamad BARBARA\at Institut de Math\'ematiques de Bourgogne(IMB)- UMR 5584 CNRS\\
Universit\'e de Bourgogne\\
9 avenue Alain Savary\\
BP 47870, 21078 Dijon cedex, France\\
\email{abdessamad.barbara@u-bourgogne.fr}
}

\maketitle

\begin{abstract}
In this paper we address a practical aspect of differential barrier penalty functions in linear programming. In this respect we propose an affine scaling interior point algorithm based on a large classe of differential barrier functions. The comparison of the algorithm with a vesion of the classical affine scaling algorithm  shows that the algorithm is robust and efficient. We thus show that differential barrier functions  open up new perspectives in linear optimization.

\end{abstract}




\section{Introduction}\label{int}

In this paper we present an algorithm based on a family of penalty functions introduced in \cite{ba}. Contrary to the classical logarithmic barrier function, these functions are not necessarily  barriers, since they can be well defined on the positive orthant including its boundary. But they are differentially barriers. In fact, these functions generalize the notion of barrier functions since (Proposition 17 of \cite{ba}) a barrier function is in particular a differential barrier one.
We recall that (Definition 1 of \cite{ba}) a function $F$ is said to be a differential barrier on the positive orthant ${\cal P}=[0,+\infty)^n$ if $F$ is differentiable on $(0,+\infty )^n $ and $\limsup\limits_{{x\to x'}\atop {x>0}}||\nabla F (x)||=+\infty$, for every $x'$ being on the boundary of ${\cal P}$. So $\nabla F$ plays the role of a barrier.
Also, the fact that  a method based on the minimization of a penalty function is of interior points type is closely related to the following
property.

\begin{proposition}\label{intdifbar}(Proposition 18  of \cite{ba})

Let $F$ be a convex, lower semi-continuous and differential barrier function on
${\cal P}$. Then every optimal solution ${\overline x}$ to the problem
$$
\min \{F(x)\ :\  Ax=b,\ x\geq 0\}$$ is an interior point of the
positive orthant.
\end{proposition}

Through the example of the concave gauge functions\footnotemark[1]\footnotetext[1]{A background on concave gauge functions is given in \cite{ba} and a complete description is done in \cite{bacr}} we will consider, we will show the important role that penalty functions of the differential barrier type can play as an alternative to the classical logarithmic barrier function.
In this respect we consider the familiy of differential barrier functions builded from the following concave gauge functions:
$$\xi_r: x\mapsto\left\{\begin{array}{ll}\big(\sum x_i^r\big)^{1\over r}&\mbox{if }
x\in[0,+\infty)^n,\cr
-\infty&\mbox{elsewhere,}
\end{array}\right.
$$
where $r$ is taken arbitrary in $(0,1)$. To be more precise,
let us consider the linear program given by

$
\hfill\min \{\left\langle c,x\right\rangle :Ax=b,\ x\in [0,+\infty)^n\
\}\hfill (LP)
$

where $A$ is an $m\times n$ matrix of rank $m$, $c$, $x\in \mathbb{R}^n$ and $b\in \mathbb{R}^m$. By definition of a concave gauge function the positive cone can be expressed as
$$[0,+\infty)^n=\left\{x\in\mathbb{R}^n:\ \xi_r(x)\geq0\right\}.$$
Hence the original linear program can be equivalently rewritten as
$$\min \{\left\langle c,x\right\rangle :Ax=b,\ \xi_r(x)\geq 0
\}.$$
Applying the approach developed in \cite{ba}, we propose to penalize the constraint $\xi_r(x)\geq 0$ by the functions
$$g_{r}:x\mapsto\left\{\begin{array}{ll} -{1\over r}\big(\xi_r(x)\big)^r&\mbox{if }x\in[0,+\infty)^n,\cr
+\infty&\mbox{elsewhere,}
\end{array}\right.$$
So the nonlinear optimization problem approximating the linear program\footnotemark[2]\footnotetext[2]{We recall that the idea to approximate a linear program by a nonlinear optimization problem is du originally to Courant \cite{courant} in 1941 with a penalty function of exterior type and later to Frisch \cite{Frisch}, in 1955, when he introduced the logarithmic barrier function which is an interior penalty one. We recall also that the notion of interior penalty operators were introduced by Auslender \cite{Auslender} in 1976 to generalize the concept of barrier functions.} is as follow

$\hfill\min\{F_{r,\mu}(x):\ Ax=b\}\hfill (P_{\mu,r})_{\mu>0}$

where
$$F_{r,\mu}(x)=\left\{\begin{array}{ll}\left\langle c,x\right\rangle+\mu g_{r}(x)&\mbox{if }x\in[0,+\infty)^n,\cr
+\infty&\mbox{elsewhere.}\end{array}\right.$$
It is easy to see that $F_{r,\mu}$ is a differential barrier function and then $\Big($Proposition \ref{intdifbar}$\Big)$ the optimal solution of $(P_{\mu,r})$ belongs to the interior of the positive orthant.

The algorithm we build, called galpv4, is of primal type and uses an affine scaling
approach\footnotemark[3]\footnotetext[3]{An affine scaling algorithm was originally proposed by Dikin in 1967 \cite{dikin}. It was rediscovered by several researchers such as barnes \cite{barnes} and Vanderbei et al \cite{vanderbeial} after Karmarkar \cite{kar} proposed his famous projective scaling algorithm.}.  It consists of two combined phases. The first one improves the feasibility of the current point and the second brings the point closer to an optimal solution. At each iteration, this requires the computation of two directions. The direction $d^k$, bringing a current point $x^k$ closer to the optimal solutions' set is obtained as follows. We compute at $x^k$ the Newton direction $d^k(\mu)$ for the problem $(P_{\mu,r})$. Vector $d^k$ is then the part of the expression of $d^k(\mu)$, independent of parameter $\mu$ and satisfies $\delta_\mu d^k(\mu)= d^k+O(\mu)$, where $\delta_\mu$ is a positive real function of $\mu$. That is $d^k=\lim\limits_{\mu\downarrow0}\delta_\mu d^k(\mu)$. The direction $d'^k$ that improves the feasibility of the current point is obtained by using the same process with the linear program
$$\min\{\lambda:\ Ax+\lambda\left(b-Ax^k\right)=b,\ x\in[0,+\infty)^n\mbox{ and }\lambda\in[0,+\infty)\}.$$
We show that the sequence $(x^k)$ converges. Its limit is an interior point of the optimal solutions' face of the linear program when $\beta\in\left(0,\displaystyle{2\over 3}\right)$, where $\beta$ is the factor of the maximal step size with respect to $x^k$ and $d^k$. Moreover, by calculating $d^k$, the algorithm generates a sequence $(y^k,s^k)$  that converges to the $\xi^\oplus_r$-analytic center of the dual optimal solutions' face, where $\xi_r^\oplus$ is the polar concave gauge function (see \cite{bacr}) of $\xi_r$.

In Proposition \ref{limitdexir} of Section \ref{algo}, we show that galpv4 includs the classical affine scaling approach by setting $r=0$. In this respect, we compare the algorithm's performances between different values of $r\in[0,1)$ through numerical experiments using the familiar netlib test set \cite{netlib}.

The paper is organized as follows. In section \ref{algo} we present the algorithm, the computation of an affine scaling  direction and how to find approximately a relative interior feasible solution. Section \ref{pstop} deals with the convergence results and a stopping criteria, followed by numerical results and comments in Section \ref{numerical}. Finally, we close the paper by some concluding remarks in Section \ref{remarks}.

\section{Presentation of a primal affine scaling method}\label{algo}

In order to take account of possible bound constraints, we consider in all the following, the linear program 

$
\hfill\min \{\left\langle c,x\right\rangle :Ax=b,\ x\geq 0,\ x_i\leq u_i\
i\in {\cal I}\}  \hfill(LPB)
$

where ${\cal I}$ is a subset to $\{1,2,\cdots,n\}$ and $u\in\mathbb{R}^n$ is given such that $u_i>0$ if $i\in {\cal I}$ and $u_i=+\infty$ if not. It is easy to see that the dual problem of $(LPB)$ is

$\hfill\max\{\langle b,y\rangle-\langle u_{\cal I},w_{\cal I}\rangle:\ A^ty+s-w=c,\ w_{\overline{\cal I}}=0, s,w\geq 0\}\hfill (LDB)$

where $\overline{\cal I}=\{1,2,\cdots,n\}-{\cal I}$. Moreover if $(x,y,s,w)$ is a primal-dual optimal solution then it is easy to see from the KKT optimality conditions that
\begin{eqnarray}\label{kktldb}
s=c-A^ty+w,\ w_{\cal I}=-U_{\cal I}^{-1}X_{\cal I}(c-A^ty)_{\cal I}\mbox{ and }w_{\overline{\cal I}}=0
\end{eqnarray}
where $U_{\cal I}=diag\left(u_{\cal I}\right)$ and $X_{\cal I}=diag\left(x_{\cal I}\right)$. We assume that there exists a relative interior feasible solution to $(LPB)$ and that the minimum is finite. Hence the optimal solutions' set of $(LPB)$ and of $(LDB)$ are non empty, and there is no duality gap. 
%
%
Moreover the set of optimal solutions to $(LDB)$ is compact.
Now taking account of the slack variables $u_i-x_i$, we adapt the definition of $\xi_r$ as
$$\xi_r:x\mapsto \left\{\begin{array}{ll}\left(\sum x_i^r+\sum\limits_{i\in\cal U}(u_i-x_i)^r\right)^{1\over r}&\mbox{if }x\geq0\mbox{ and }u_i-x_i\geq0,\ \forall i\in{\cal I},\cr
-\infty&\mbox{elsewhere.}\end{array}\right.$$


The following algorithm, called {\it galpv4}, uses an approach based on a version of the classical affine scaling algorithm presented in \cite{barnes,vanderbeial,saigal}. 
$\underline {Algorithm}$ galpv4
\\
$\underline {Initialization}$

Construct a starting point $x$ as described just bellow and choose $r\in[0,1)$.

Compute $y$, $w$, $s$ according to (\ref{dualy}), (\ref{dualw}) and (\ref{duals}) respectively. 

Compute the expected relative duality gap $Rgap$ according to (\ref{rgap})

Set the feasibility measure $Rf\leftarrow\displaystyle{\|Ax-b\|_\infty\over \|b\|_\infty+1}$

Choose $\epsilon>0$ a stopping rule parameter.\\

{\bf While} $\min(Rf,Rgap)>\epsilon$ {\bf do}

\hskip 1cm Compute  $d_x$ the feasible direction according to (\ref{dx})

\hskip 1cm Compute $d$ the descent direction according to (\ref{d})

\hskip 1cm Set $t_{max}\leftarrow\min\left(\min\limits_{d_{x_i}<0}-\displaystyle {x_i\over d_{x_i}},
\min\limits_{{d_{x_i}>0},{i\in{\cal I}}}-\displaystyle
{u_i-x_i\over d_{x_i}},1\right)$

\hskip 1cm If $Rf> \epsilon$  then $t\leftarrow0.95t_{max}$ else set $t\leftarrow0.65t_{max}$

\hskip 1cm Set $x\leftarrow x+td_x$

\hskip 1cm Set $t_{max}\leftarrow\min\left(\min\limits_{d_i<0}-\displaystyle {x_i\over d_i},
\min\limits_{{d_i>0},{i\in{\cal I}}}-\displaystyle {u_i-x_i\over d_i}\right)$

\hskip 1cm If $Rf> \epsilon$  then $t\leftarrow0.65t_{max}$ else set $t\leftarrow0.95t_{max}$

\hskip 1cm Update $x\leftarrow x+td$

\hskip 1cm Update $y$, $w$, $s$ according to (\ref{dualy}), (\ref{dualw}) and (\ref{duals}) respectively. 

\hskip 1cm Update the expected relative duality gap $Rgap$ according to (\ref{rgap})

\hskip 1cm Update $Rf\leftarrow\displaystyle{\|Ax-b\|_\infty\over \|b\|_\infty+1}$

{\bf End while}

Let us describe now how to construct, empirically, a starting point. In fact we construct two starting points $x^1$ and $x^2$.  
The first one is defined as follow. $\mbox{For }j=1,2,..,n$, $x_j^1=
\min\left(\displaystyle {n\over \left\|A_{j.}\right\|},0.9u_j\right)$ if $c_j<0$ and $x_j^1=
\min\left(\displaystyle {n\over \left\|A_{j.}\right\|},0.1u_j\right)$ otherwise, where  $A_{j.}$ is the $j^{th}$ column of matrix $A$. The second one is defined as in the routine $pcinit.f$ of the software HOPDM of Gondzio \cite{hopdm}. Our starting point $x$ is chosing as follow. If $\min x^2_i>\min\limits_j x_j^1$ or $\min\limits_j x_j^1<1$ then we set $x=x^2$ else we set $x=x^1$.

Note that the algorithm can be extended to the case $r=0$. It is justified by the following proposition.

\begin{proposition}\label{limitdexir}

Set $n_{\cal I}=card({\cal I})$ and define  $\xi_0$ as
$$\xi_0(x)=\displaystyle\left\{\begin{array}{ll}\left(\prod\limits_{i\in\{1,\cdots,n\}} x_i\prod\limits_{i\in{\cal I}} (u_i-x_i)\right)^{1\over n+n_{\cal I}}&\mbox{if }x\in[0,+\infty)^n\cr
-\infty&\mbox{elsewhere}.\end{array}\right.
$$
For any $r\in (-\infty,0)\cup (0,1)$ we set ${\tilde \xi_r}={1\over (n+n_{\cal I})^{1\over r}}\xi_r$. Then for every $x\in (0,+\infty )^n$,

{\bf (i)} $\lim\limits_{r\uparrow 0}{\tilde \xi_r}(x)=
            \lim\limits_{r\downarrow 0}{\tilde \xi_r}(x)={1\over n+n_{\cal I}}\xi _0(x)$,

{\bf (ii)} $\lim\limits_{r\uparrow 0}\nabla{\tilde \xi_r}(x)=
            \lim\limits_{r\downarrow 0}\nabla{\tilde \xi_r}(x)={1\over n+n_{\cal I}}\nabla\xi
_0(x)$,

{\bf (iii)} $\lim\limits_{r\uparrow 0}\nabla ^2{\tilde \xi_r}(x)=
            \lim\limits_{r\downarrow 0}\nabla ^2{\tilde \xi_r}(x)={1\over n+n_{\cal I}}\nabla
^2\xi _0(x)$.
\end{proposition}

{\it Proof.} {\bf (i).} Without loss of generality we can assume that ${\cal I}=\emptyset$. Let $x\in (0,+\infty )^n$ and set
$\psi _1(r)=\ln({1\over n}\Sigma x_i^r)\hbox{ and }\psi _2(r)=r.$
We have
$\lim\limits_{r\to 0}\psi _1(r)=\lim\limits_{r\to 0}\psi
_2(r)=0\hbox{ and }\psi '_2(r)=1\not =0.$
 Then by the classical H\^opital theorem
$\lim\limits_{r\to 0}\ln{\tilde \xi_r}(x)=\displaystyle
\lim\limits_{r\to 0}{\psi _1(r)\over \psi _2(r)}
=\displaystyle\lim\limits_{r\to 0}{\psi '_1(r)\over \psi '_2(r)} ,$
but
$\psi '_1(r)=\displaystyle{
\sum\limits_{i=1}^nx_i^r\ln x_i\over
\sum\limits_{i=1}^nx_i^r}. $
The result follows.

{\bf(ii)} and {\bf(iii).} Using { (i)} and the expressions of $\nabla {\tilde \xi_r}$ and $\nabla^2{\tilde \xi_r}$, it is easy to see that
$\lim\limits_{r\to 0}\nabla{\tilde \xi _r}(x)=
\nabla{1\over n}\xi _0(x)\hbox{ and }
\lim\limits_{r\to 0}\nabla ^2{\tilde \xi _r}(x)=
\nabla ^2{1\over n}\xi _0(x).
$
\hfill\qed

\subsection{\bf Finding a descent direction}\label{N}

Let $x$ be a relative interior feasible point to $(LPB)$, $\mu>0$ and $r\in(0,1)$. The Newton direction at $x$ to the penalized problem  $(P_{\mu,r})$ is obtained by solving the minimization problem
$$\min\left\{ {\displaystyle 1\over 2}\left\langle \nabla ^2F_{r,\mu} (x) d,d\right\rangle +
\left\langle \nabla F_{r,\mu} (x), d\right\rangle\ :\ Ad=0\right\} .$$ Using
the KKT optimality conditions, the problem amounts to finding $d(\mu
)\in \mathbb{R} ^n$ and $y\in \mathbb{R}^m$ solutions to the system of linear
equations

\begin{equation}\label{enewt}\left\{ \begin{array}{*{1}c@{\;=\;}c}
\nabla ^2F_{r,\mu} (x) d(\mu )+\nabla F_{r,\mu} (x)+A^ty   & 0 \cr
                             \hfill       Ad(\mu )  &  0.
\end{array}\right.
\end{equation}

We have
$
\nabla F_{r,\mu} (x)=c - \mu Ge\mbox{ and }\nabla ^2F_{r,\mu} (x) =  \mu (1-r)H
$
where $e^t=(1,1,..,1)\in\mathbb{R}^n$,
$G$ and $H$ are diagonal matrices defined respectively by
\begin{equation}\label{G}
G_{ii}=\left\{ \begin{array}{*{1}c@{\; \;}c}
x_i^{r-1} & \mbox{ if }i\notin {\cal I},\\
x_i^{r-1}-(u_i-x_i)^{r-1} & \mbox{ otherwise }
\end{array}\right.
\end{equation}
and
\begin{equation}\label{H}
H_{ii}=\left\{ \begin{array}{*{1}c@{\; \;}c}
x_i^{r-2} & \mbox{ if }i\notin {\cal I},\\
x_i^{r-2}+(u_i-x_i)^{r-2} & \mbox{ otherwise. }
\end{array}\right.
\end{equation}
Then setting
\begin{equation}\label{P}
P=I-H^{-{1\over 2}}A^t\big(AH^{-1}A^t\big)^{-1}AH^{-{1\over 2}},
\end{equation}
the projection matrix on the kernel of $AH^{-{1\over 2}}$, system
(\ref{enewt}) reduces to
$PH^{-{1\over 2}}\nabla F_{r,\mu} (x) 
+ \mu (1-r)H^{1\over 2}d(\mu ) =0
$
and then
$\mu (1-r)d(\mu ) =
-H^{-{1\over 2}}PH^{-{1\over 2}}\nabla F_{r,\mu} (x)=-H^{-{1\over 2}}PH^{-{1\over 2}}c+\mu H^{-{1\over 2}}PH^{-{1\over 2}}Ge.$
Since $\mu(1-r)>0$ we can so take as an affine scaling direction at $x$ to the linear program vector $d$ given by

\begin{equation}\label{d}
d=\lim\limits_{\mu\downarrow0}\mu (1-r)d(\mu ) =-H^{-{1\over 2}}PH^{-{1\over 2}}c
\end{equation}

Observe that since $\langle c,d\rangle=-\|PH^{-{1\over 2}}c\|_2^2<0$, $d$ is a descent direction to the linear program at every point to $\mathbb{R}^n$.

{\bf Remark:} To improve the quality of the direction $d$, in order to maintain a good feasibility to the current point, we can compute in addition the direction $H^{-{1\over 2}}PH^{{1\over 2}}d$ which can be used instead of $d$. Which in fact amounts to projecting a second time the direction $H^{{1\over 2}}d$ onto the vector subspace $\ker AH^{-{1\over 2}}$. Of course, theoretically the two directions are identical, but numerically there is a significant difference. However the computation of $H^{-{1\over 2}}PH^{{1\over 2}}d$ has some extra cost in number of operations. Therefore we use the technic only when the relative duality gap is less than 0.001 or the current number of iterations exceeds 20.


\subsection{\bf Finding a feasible solution}\label{feasiblesolution}

It is well known that an approximate relative interior feasible solution to $(LPB)$ can be obtained by solving a linear problem
of the form

$\hfill\min\left\{\lambda  :Ax+\lambda (b-Ax^0)=b,\ x\geq 0,\ x_i\leq u_i,\ \forall i\in{\cal I},\ \lambda\geq0
\right\} , \hfill(FLP )
$

where $x^0$ is a point arbitrarily chosen in $(0,+\infty)^n$.
Write $(FLP)$ as
$$\min\left\{\left\langle {\tilde c},{\tilde x}\right\rangle :\ {\tilde A}{\tilde x}=b,\ {\tilde x}\geq 0,\ {\tilde x}_i\leq u_i,\ \forall
i\in{\cal I}\right\} ,
$$
where ${\tilde x}=\left(\begin{array}{l}x\cr\lambda\end{array} \right),\ {\tilde c}=\left(\begin{array}{l}0_{\mathbb{R}^n}\cr 1\end{array}\right)\mbox{ and
} {\tilde A}=\left( \matrix{ A & b-Ax^0\cr }\right) .
$
Then using (\ref{d}), the affine scaling direction ${\tilde d}$ with respect to ${\tilde x}$ is given by
${\tilde d} = -{\tilde H}^{-{1\over 2}}{\tilde P}{\tilde H}^{-{1\over 2}}
{\tilde c}
$
where ${\tilde H}=\left( \matrix{ H & 0\cr 0 & \lambda^{r-2}\cr }\right)
\mbox{ and  }{\tilde P}=I-{\tilde H}^{-{1\over 2}}{\tilde A}^t\left({\tilde
A}{\tilde H}^{-1} {\tilde A}^t\right)^{-1}{\tilde A}{\tilde H}^{-{1\over
2}}. $

But the matrix ${\tilde A}{\tilde H}^{-1}{\tilde A}^t$ will be
generally dense when there is one dense column in ${\tilde A}$. Column
$b-Ax^0$, in most cases, is dense. So for large-scale applications,
we split such column from the others. We proceed as follow. Set
$v=b-Ax^0$. Then ${\tilde A}{\tilde H}^{-1}{\tilde
A}^t=AH^{-1}A^t+\lambda^{2-r}vv^t.$ Using the Sherman-Morrison
formula we have $ \left({\tilde A}{\tilde H}^{-1}{\tilde
A}^t\right)^{-1}=\left(AH^{-1}A^t\right)^{-1}+\delta \lambda^{2-r}ww^t, $ where
$w=\left(AH^{-1}A^t\right)^{-1}v$
and $ \delta=\displaystyle{-1\over 1+\lambda^{2-r}\left\langle
w,v\right\rangle }. $ So

$$\begin{array}{ll} {\tilde P}&=\left( \matrix {
P-\delta\lambda^{2-r}H^{-1\over 2}A^tww^tAH^{-1\over
2}&\delta\lambda^{1-{r\over 2}} H^{-1\over 2}A^tw\cr&\cr
\delta\lambda^{1-{r\over 2}}w^tAH^{-1\over 2}& -\delta \cr }
        \right)\cr
&\cr&=\left( \matrix {
P & 0\cr&\cr
0 & 0 \cr }
        \right)
    -\delta \left( \matrix {
\lambda^{1-{r\over 2}}H^{-1\over 2}A^tw\cr&\cr
-1 \cr }
        \right)
\left( \matrix {
\lambda^{1-{r\over 2}}w^tAH^{-1\over 2} && -1\cr }
        \right) .
\end{array}$$
It follows that
$
{\tilde d}=-\delta\lambda^{2-r}\left( \matrix{H^{-1}A^tw
\cr\cr&\cr
                             -1    \cr}\right).
$
\noindent Since $-\delta\lambda^{2-r}>0$ the search directions with respect to $x$ and $\lambda$ can be expressed respectively as
 $d_x=-H^{-1}A^t\left(AH^{-1}A^t\right)^{-1}(b-Ax^0)\mbox{ and }d_\lambda = -1.$
\noindent But numerical experiments show that as iterations go, the
constraint $Ax+\lambda \left(b-Ax^0\right)=b$ is less and less satisfied. 
This is due to the rounding off errors generated by the projection onto $\ker\tilde{A}$ at each iteration and thus creating a snowball effect. To work around this problem, we
proceed as follows: Let $x^k$ be a current point in the
feasibility searching phase. Then $\left(\begin{array}{l}x^k\cr1\end{array}\right)$ is a feasible point of
problem
$$\min\left\{ \lambda :\ Ax+\lambda \left(b-Ax^k\right)=b,\ x\geq 0,\ x_i\leq u_i,\ \forall i\in {\cal I},\  \lambda\geq 0\right\} .
$$
\noindent In this case the search direction with respect to $x^k$ is 
\begin{eqnarray}\label{dx}
d_x^k=-H^{-1}A^t\left(AH^{-1}A^t\right)^{-1}(b-Ax^k)
\end{eqnarray} 
It follows that the point
$\left(\begin{array}{l}x^{k+1}\cr\lambda^{k+1}\end{array}\right)=\left(\begin{array}{l}x^k\cr1\end{array}\right)+t^k\left(\begin{array}{l}{ d_x^k}\cr-1\end{array}\right)$ for a step size
$t^k$ suitably chosen, does not suffer from the snowball effect mentioned above.

{\bf Remark:} To compute $\left(AH^{-1}A^t\right)^{-1}\left(b-Ax\right)$ and $PH^{-{1\over
2}}c$ we solve for $w$ and $\Delta$ by Cholesky factorization the linear systems
$AH^{-1}A^tw=b-Ax$ and $\left(AH^{-1}A^t\right)\Delta=H^{-{1\over 2}}c$ and then we compute
$PH^{-{1\over 2}}c = H^{-{1\over 2}}c - H^{-{1\over 2}}A^t\Delta$.

\section{Convergence, dual solution and stopping criteria}\label{pstop}

Without loss of generality we can assume in this section that ${\cal I}=\emptyset$. In this case
$H=X^{r-2}\mbox,\ P=I-X^{1-{r\over
2}}A^t(AX^{2-r}A^t)^{-1}AX^{1-{r\over 2}},
$
 where $X=diag(x)$ and $x\in(0,+\infty)^n$. 
 To simplify we assume that the starting point $x^0$ is a relative interior feasible solution to the linear program. So we consider $(x^k)_{k\in\mathbb{N}}$ the sequence defined by
$x^{k+1}=x^k+\beta t^k_{max}d^k,$ where $\beta\in(0,1)$ and $t^k_{max}$ is the maximum step length with
respect to $x^k$ and $d^k=-X_k^{1-{r\over 2}}PX_k^{1-{r\over 2}}c=-X_k^{2-r}c+X_k^{2-r}A^t(AX_k^{2-r}A^t)^{-1}AX_k^{2-r}c$. Set $y^k=\left(AX_k^{2-r}A^t\right)^{-1}AX_k^{2-r}c\mbox{ and }s^k=c-A^ty^k.$ Here is our main result.

\begin{theorem}\label{convergenceys}
Assume that $\displaystyle0<\beta<{2\over 3}$. Then $\left(x^k,y^k,s^k\right)_{k\in\mathbb{N}}$ converges to $\left(\overline{x},{\overline y},{\overline s}\right)$, where $\left({\overline y},{\overline s}\right)$ is the $\xi_t$-analytic center to the dual optimal face of the linear program, $t$ is such that $\displaystyle{1\over t}+{1\over r}=1$ and ${\overline x}$ belongs to  the relative interior of the primal optimal face of the linear program. 
\end{theorem}

Before giving the proof of the theorem, we first establish some preliminary results.

\begin{proposition}\label{series}

$\sum\limits_{k=0}\limits^{\infty} \beta t^k_{max}\left\|PX_k^{1-{r\over 2}}c\right\|_2^2$ is a converging series.

\end{proposition}

{\it Proof.} 
We have $\left\langle c,x^{k+1}\right\rangle -\left\langle c,x^k\right\rangle= \beta t^k_{max}\left\langle
c, d^k\right\rangle=- \beta t^k_{max}\left\|PX_k^{1-{r\over 2}}c\right\|_2^2.$
The sequence $\left(\left\langle c,x^k\right\rangle\right)_{k\in\mathbb{N}}$
is then decreasing. Since we assumed that the optimal value of the linear program is finite, the sequence is bounded and then converges. Set ${\overline c}$ its limit.
Then we have $\sum\limits_{k=0}\limits^{\infty}  \beta t^k_{max}\left\|PX_k^{1-{r\over 2}}c\right\|_2^2=\left\langle c,x^0\right\rangle -{\overline c}<+\infty.$ The result then follows.
\hfill\qed

Now let us recall an important result. It was proved by Monteiro et al. \cite{monteiroetal}, Saigal \cite{saigal}, Tseng and Luo \cite{tsenluo} and Tsuchiya \cite{tsuchia}.

\begin{theorem}\label{résultatimportant}
There exists a constant $L(A,c)>0$ such that every optimal solution ${\overline w}$ to the following ellipsoidal problem

$\hfill\max\left\{ \left\langle c,w\right\rangle:\ Aw=0,\ \left\|X^{-1}w\right\|_2^2\leq 1\right\}\hfill(EP)$

satisfies $\left\|{\overline w}\right\|_2\leq L(A,c)\left\langle c,{\overline w}\right\rangle$.
\end{theorem}

\begin{corollary}\label{résultatimportantc}
Let $x\in(0,+\infty)^n$. Then $X^{1-{r\over 2}}PX^{1-{r\over 2}}c$ satisfies
$$\left\|X^{1-{r\over 2}}PX^{1-{r\over 2}}c\right\|_2\leq L(A,c)\left\langle c,X^{1-{r\over 2}}PX^{1-{r\over 2}}c\right\rangle.$$
\end{corollary}

{\it Proof.} First observe that $\displaystyle X^{1-{r\over 2}}PX^{1-{r\over 2}}c\over \displaystyle \left\|PX^{1-{r\over 2}}c\right\|_2$ can be viewed as the optimal solution to the following ellipsoidal problem

$\hfill\max\left\{ \left\langle c,w\right\rangle:\ Aw=0,\ \left\|X^{{r\over 2}-1}w\right\|_2^2\leq 1\right\}\hfill(EP_r)$

Hence using Theorem \ref{résultatimportant} by considering ${\tilde X}=X^{1-{r\over 2}}$ instead of $X$, the result follows.
\hfill\qed

\begin{proposition}\label{convergencex}

$(x^k)_{k\in\mathbb{N}}$ is a converging sequence, say to ${\overline x}$. Furthermore, for each $k\in\mathbb{N}$, $\left\|x^k-{\overline x}\right\|_2\leq h\left(\left\langle c,x^k\right\rangle -{\overline c}\right),$ where $\displaystyle h={1\over L(A,c)}$.
\end{proposition}

{\it Proof.} By Corollary \ref{résultatimportantc} we have
$\left\langle c,x^{k}\right\rangle -\left\langle c, x^{k+1}\right\rangle=- \beta t^k_{max}\left\langle c,d^k\right\rangle\geq L(A,c)\left\| \beta t^k_{max}d^k\right\|_2=L(A,c)\left\|x^{k+1}-x^k\right\|_2.$
It follows that
$+\infty>\left\langle c,x^0\right\rangle-{\overline c}=\sum\limits_{0\leq k<+\infty}\left\langle c,x^k
-x^{k+1}\right\rangle\geq L(A,c)\sum\limits_{0\leq k<+\infty}\left\|x^{k+1}-x^k\right\|_2.$ Thus $(x^k)_{k\in\mathbb{N}}$ converges.
Now using again Corollary \ref{résultatimportantc} we have
$\displaystyle\left\langle c, x^k\right\rangle-{\overline c}=\sum\limits_{j=0}^{\atop \infty}\left\langle c,x^{k+j}
-x^{k+j+1}\right\rangle\geq {1\over h}\sum\limits_{j=0}^{\atop \infty}\left\|x^{k+j}-x^{k+j+1}\right\|_2
\geq \displaystyle{1\over h}\left\|\sum\limits_{j=0}^{\atop \infty}x^{k+j}-x^{k+j+1}\right\|_2={1\over h}\left\|x^k-{\overline x}\right\|_2.$
The result follows.
\hfill\qed

Now we recall the next theorem proved by Dikin \cite{dikinbis}. A proof can also be found in \cite{saigal,todd,vandlaga,vavasis}.

\begin{theorem}\label{thdikin}
For every $x>0$ and for every $p\in\mathbb{R}^n$, we have 
$$\left\|\left(AX^2A^t\right)^{-1}AX^2p\right\|_2\leq q(A)\left\|p\right\|_2,$$ 
where $q(A)$ is a constant only function of $A$.
\end{theorem}


\begin{proposition}\label{yksk}

The sequences $(y^k)$ and $(s^k)$ are bounded.
\end{proposition}

{\it Proof}
According to Theorem \ref{thdikin}, for every $x>0$ and for every $p\in\mathbb{R}^n$, we have $\|y^k\|_2=\left\|\left(AX_k^{2-r}A^t\right)^{-1}AX_k^{2-r}c\right\|_2=\left\|\left(A\left(X_k^{1-{r\over 2}}\right)^2A^t\right)^{-1}A\left(X_k^{1-{r\over 2}}\right)^2c\right\|_2\leq q(A)\left\|c\right\|_2$ and then $\|s^k\|_2=\|c-A^ty^k\|_2\leq(1+q(A)\|A\|_2)\|c\|_2$. The result then follows.\hfill\qed

Let us now consider the following notation. Given $x\in(0,+\infty)^n$ and $s\in\mathbb{R}^n$ we set $I_r(x,s)=\{i:\ x_i^{1-r}|s_i|=\|X^{1-r}s\|_\infty\}$.

\begin{lemma}\label{prelim1.1}

Let $(x,s)\in(0,+\infty)^n\times\mathbb{R}^n$ be such that $Xs\not=0$. One has for every $(r,r')\in[0,1]^2$, if $r<r'$ then $x_{i_r}\geq x_{i_{r'}}$ and $|s_{i_r}|\leq |s_{i_{r'}}|$, $\forall(i_r,i_{r'})\in I_r(x,s)\times I_{r'}(x,s)$.
\end{lemma}

{\it Proof.} We have
$$0<x_{i_{r'}}^{1-r}|s_{i_{r'}}|\leq\|X^{1-r}s\|_\infty=x_{i_{r}}^{1-r}|s_{i_{r}}|\ \ \ \ \ \ \ (1)$$
and
$$0<x_{i_{r}}^{1-r'}|s_{i_{r}}|\leq\|X^{1-r'}s\|_\infty=x_{i_{r'}}^{1-r'}|s_{i_{r'}}|\ \ \ \ \ (2)$$
Multiplying side by side $(1)$ and $(2)$ one has
$0<x_{i_{r'}}^{1-r}x_{i_{r}}^{1-r'}|s_{i_{r}}||s_{i_{r'}}|\leq x_{i_{r}}^{1-r}x_{i_{r'}}^{1-r'}|s_{i_{r'}}||s_{i_{r}}|$.
That is $x_{i_{r'}}^{r'-r}\leq x_{i_{r}}^{r'-r}$ and then $x_{i_{r'}}\leq x_{i_{r}}$. Now using $(2)$ one has
$0<x_{i_{r'}}^{1-r'}|s_{i_{r}}|\leq x_{i_{r}}^{1-r'}|s_{i_{r}}|\leq \|X^{1-r'}s\|_\infty=x_{i_{r'}}^{1-r'}|s_{i_{r'}}|$. The result then follows.\hfill\qed

Define $I=\{i:\ {\overline x}_i=0\},\ J=\{i:\ {\overline x}_i>0\}\mbox{ and }n_I=card(I).$

\begin{lemma}\label{prelim2.1}
There is ${\tilde h}> 0$ such that $\|x^k_J-{\overline x}_J\|_2\leq{\tilde h}\|x^k_I\|_2,\ \forall k\in\mathbb{N}.$

\end{lemma}

{\it Proof.}
Let $({\overline y},{\overline s})$ be an accumulation point of $(y^k,s^k)$. The existence of $({\overline y},{\overline s})$ is ensured by Proposition \ref{yksk}. Using Proposition \ref{convergencex} we have
$\|x^k-{\overline x}\|_2^2=\|x^k_I\|_2^2+\|x^k_J-{\overline x}_J\|^2_2\leq h^2\langle c,x^k-{\overline x}\rangle^2=h^2\langle {\overline s},x^k-{\overline x}\rangle^2=h^2\langle {\overline s}_I,x^k_I\rangle^2\leq h^2\|{\overline s}_I\|^2_2\|x^k_I\|_2^2 $
and then $\|x^k_J-{\overline x}_J\|_2^2\leq(h^2\|{\overline s}_I\|^2_2-1)\|x^k_I\|_2^2.$ Thus $h^2\|{\overline s}_I\|^2_2-1$ is necessarily nonnegative. The result then follows by setting $\tilde{h}=\sqrt{h^2\|{\overline s}_I\|^2_2-1}$.\hfill\qed

{\bf N.B:} The fact that $h=\displaystyle{1\over L(A,c)}$, $h^2\|{\overline s}_I\|^2_2-1\geq 0$ also means that $\|\overline{s}\|_2\geq L(A,c)$, for every $\overline{s}$ be an accumulation point of $(s^k)$.

In all the following we set
$g=\sup\limits_{k\in\mathbb{N}}\|s^k\|_\infty$, ${\overline M}=\sup\limits_{k} \|x^k\|_\infty<+\infty$ and $\underline{M}=\inf\limits_{k}\min\limits_{j\in J}x_j^k$. Note that since $\lim\limits_{k\uparrow+\infty}x_J^k={\overline x_J}>0$, $\underline{M}>0$ and that according to Proposition \ref{yksk} $g<+\infty$. 

\begin{proposition}\label{convergencecx}
Let $\beta\in(0,1)$. Then there exists $K\in\mathbb{N}$ such that

{\bf i)} $\forall k\geq K$, $\forall r\in(0,1)$, $\forall i_r\in I_r(x^k,s^k)$, $x_{i_r}^k=O(\|x_I^k\|_\infty)$ and $s_{i_r}^k=O(1)$. Furthermore $\|X_k^{1-r}s^k\|_\infty=\|X_{k,I}^{1-r}s_I^k\|_\infty=O\left(\|x_I\|^{1-r}_\infty\right)$ and there exists a constant $\hat{C}$ such that $\left\|s_J^k\right\|_2\leq \hat{C}\|x_I\|_\infty^{2-r}$.

{\bf ii)} $\forall k\geq K$, $\left\langle c,x^{k+1}\right\rangle -{\overline c}\leq{\overline L}(\langle c,x^k\rangle-{\overline c})$, 
where $\displaystyle{\overline L}=1-\beta{L(A,c)^{6-r\over 2-r}\over 2^{3-r\over 2-r}g^{6-2r\over 2-r}n_I^{{7\over 2}-{r(1-r)\over 2(2-r)}}}$.

{\bf iii)} $\sum\limits_{k=0}^{\atop \infty}\|x^k_I\|_\infty^a<+\infty,\ \forall a>0.$


{\bf iv)} $\langle x^k_I,s^k_I\rangle=O(\|x_I^k\|_\infty)$.

{\bf v)} $\langle c,x^k\rangle-\overline{c}=O(\|x^k_I\|_\infty)$ and $\|x^k-\overline{x}\|_2=O(\|x^k_I\|_\infty)$.
\end{proposition}

{\it Proof.} {\bf i)} and {\bf ii)} We have
$\left\langle c,x^{k+1}\right\rangle -{\overline c}=\left\langle c,x^k\right\rangle -{\overline c}-\beta t_{max}\left\|PX_k^{1-{r\over 2}}c\right\|_2^2=\left\langle c,x^k\right\rangle -{\overline c}-\beta t_{max}\left\|X_k^{1-{r\over 2}}s^k\right\|_2^2\leq\left\langle c,x^k\right\rangle -{\overline c}-\beta t_{max}\left\|X_k^{1-{r\over 2}}s^k\right\|^2_\infty .$
Now $\displaystyle t^k_{max}=\min\left\{-{x^k_i\over d^k_i}:\ d^k_i<0\right\}\geq {1\over \|X_k^{-1}d^k\|_\infty}={1\over\|X_k^{1-r}s^k\|_\infty}$. Then $\displaystyle\left\langle c,x^{k+1}\right\rangle -{\overline c}\leq\left\langle c,x^k\right\rangle -{\overline c}-\beta {\left\|X_k^{1-{r\over 2}}s^k\right\|^2_\infty\over\left\|X_k^{1-r}s^k\right\|_\infty} .$ Let $\left(i_{r\over 2},i_r\right)\in I_{r\over 2}(x^k,s^k)\times I_r(x^k,s^k)$. From Lemma \ref{prelim1.1} one has $x^k_{i_{r\over 2}}\geq x^k_{i_r}$ and then

$\hfill\displaystyle\left\langle c,x^{k+1}\right\rangle -{\overline c}\leq\left\langle c,x^k\right\rangle -{\overline c}-\beta x^k_{i_{r\over 2}}{\left|s^k_{i_{r\over 2}}\right|^2\over \left|s^k_{i_r}\right|}\hfill (1)$

Using Proposition \ref{convergencex} , the fact that $\displaystyle X_k^{1-{r\over 2}}PX_k^{1-{r\over 2}}c\over \displaystyle \left\|PX_k^{1-{r\over 2}}c\right\|_2$ is the optimal solution to $(EP_r)$ and $\displaystyle x^k-{\overline x}\over \displaystyle \left\|X_k^{{r\over 2}-1}\left(x^k-{\overline x}\right)\right\|_2$ is a feasible solution of $(EP_r)$,
$L(A,c)\displaystyle{\|x^k-{\overline x}\|_2\over\left\|X_k^{{r\over 2}-1}(x^k-{\overline x})\right\|_2}\leq {\left\langle c,x^k-{\overline x}\right\rangle\over \left\|X_k^{{r\over 2}-1}(x^k-{\overline x})\right\|_2}\leq {\left\langle c,X_k^{1-{r\over 2}}PX_k^{1-{r\over 2}}c\right\rangle\over \left\|PX_k^{1-{r\over 2}}c\right\|_2}=\left\|PX_k^{1-{r\over 2}}c\right\|_2=\|X_k^{1-{r\over 2}}s^k\|_2.$
It follows that
$\displaystyle L(A,c)^2{\|x^k_I\|_2^2+\|x^k_J-{\overline x}_J\|_2^2\over \left\|{X_k}_I^{r\over 2}e_I\right\|_2^2+\left\|{X_k}_J^{{r\over 2}-1}(x^k_J-{\overline x}_J)\right\|_2^2}\leq \|X_k^{1-{r\over 2}}s^k\|_2^2.$ Now using Lemma \ref{prelim2.1}, the fact that $r\in(0,1)$ and the fact that $\lim\limits_{k\to\infty}x^k_I=0$, for $k$ being large enough one has
$\left\|{X_k}_J^{{r\over 2}-1}(x^k_J-{\overline x}_J)\right\|_2^2\leq \underline{M}^{r-2}\tilde{h}^2\|x_I^k\|_2^2=\underline{M}^{r-2}\tilde{h}^2\|{X_k}_I^{1-{r\over 2}}{X_k}_I^{r\over 2}e_I\|_2^2\leq\underline{M}^{r-2}\tilde{h}^2\|x_I^k\|^{2-r}_\infty\left\|{X_k}_I^{r\over 2}e_I\right\|_2^2\leq\left\|{X_k}_I^{r\over 2}e_I\right\|_2^2$, where $e_I$ is the vector of $\mathbb{R}^{n_I}$ whose components are equal to 1.
According to iii) of Proposition 4.3 in \cite{bacr} one has
$\left(\sum\limits_{i\in I}{x^k_i}^r\right)^{2\over r}n_I^{1-{2\over r}}=\left(\sum\limits_{i\in I}{({x^k_i}^2)}^{r\over 2}\right)^{2\over r}n_I^{1-{2\over r}}=\xi_{{r\over 2},I}(x_I^k)\xi_{{r\over r-2},I}(e_I)\leq\langle x_I^k,e_I\rangle= \sum\limits_{i\in I}{x^k_i}^2=\|x^k_I\|^2_2,$
where $\xi_{{r\over2},I}$ and $\xi_{{r\over r-2},I}$ are the concave gauge functions respectively defined by 
$\xi_{{r\over2},I}(z)=\left\{\begin{array}{ll}
\left(\sum\limits_{i\in I}z_i^{r\over2}\right)^{2\over r}&\mbox{if }z\in[0,+\infty)^{n_I},\cr
-\infty&\mbox{elsewhere}
\end{array}\right.$ and $\xi_{{r\over r-2},I}(z)=\left\{\begin{array}{ll}
\left(\sum\limits_{i\in I}z_i^{r\over r-2}\right)^{r-2\over r}&\mbox{if }z\in(0,+\infty)^{n_I},\cr
0&\mbox{if }z\in\partial[0,+\infty)^{n_I},\cr
-\infty&\mbox{elsewhere.}
\end{array}\right.$
Here $\partial[0,+\infty)^{n_I}$ denotes the boundary of $[0,+\infty)^{n_I}$.
That is
$\displaystyle\sum\limits_{i\in I}{x_i^k}^r=\|{X_k}_I^{r\over 2}e_I\|_2^2=\sum\limits_{i\in I}{x^k_i}^r\leq n_I^{1-{r\over 2}}\|x^k_I\|^r_2.$
Hence
$\displaystyle{{L(A,c)}^2\over 2n_I^{1-{r\over 2}}}\|x_I\|_2^{2-r}\leq L(A,c)^2{\|x^k_I\|_2^2+\|x^k_J-{\overline x}_J\|_2^2\over \left\|{X_k}_I^{r\over 2}e_I\right\|_2^2+\left\|{X_k}_J^{{r\over 2}-1}(x^k_J-{\overline x}_J)\right\|_2^2}\leq \|X_k^{1-{r\over 2}}s^k\|_2^2$
and then

$\hfill\displaystyle{{L(A,c)}^2\over n_I^{1-{r\over 2}}}\|x_I\|_2^{2-r}\leq 2n_I\|X_k^{1-{r\over 2}}s^k\|_\infty \hfill(2)$

\noindent Now using Corollary \ref{résultatimportantc} one has

\noindent$\underline{M}^{2-r}\|s^k_J\|_2\leq\underline{M}^{1-{r\over 2}}\left\|X_{k,J}^{1-{r\over 2}}s_J^k\right\|_2\leq\left(\min\limits_{j\in J}x_j^k\right)^{1-{r\over 2}}\left\|X_{k,J}^{1-{r\over 2}}s_J^k\right\|_2$

$\leq\left\|\left(X_k^{1-{r\over 2}}PX_k^{1-{r\over 2}}c\right)_J\right\|_2\leq\left\|X_k^{1-{r\over 2}}PX_k^{1-{r\over 2}}c\right\|_2\leq L(A,c)\left\langle c,X_k^{1-{r\over 2}}PX_k^{1-{r\over 2}}c\right\rangle$

$=L(A,c)\left\langle {\overline s}_I,X_{k,I}^{1-{r\over 2}}X_{k,I}^{1-{r\over 2}}s^k_I\right\rangle\leq L(A,c)\left\|{\overline s}_I\right\|_2\left\|X_{k,I}^{1-{r\over 2}}X_{k,I}^{1-{r\over 2}}s^k_I\right\|_2$

$\leq L(A,c)\sqrt{n_I}\|\overline{s}_I\|_\infty\|x_I^k\|^{1-{r\over 2}}_\infty\left\|X_{k,I}^{1-{r\over 2}}s^k_I\right\|_2\leq L(A,c)\left\|{\overline s}_I\right\|_\infty n_Ig\|x^k_I\|_\infty^{2-r}$.

\noindent So we get on the one hand $\|s^k_J\|_2\leq\hat{C}\|x^k_I\|^{2-r}$, where $\hat{C}=\displaystyle{L(A,c)\left\|{\overline s}_I\right\|_\infty n_Ig\over\underline{M}^{2-r}}$, and 

$\hfill\left\|X_{k,J}^{1-{r\over 2}}s_J^k\right\|_2\leq\displaystyle{L(A,c)\sqrt{n_I}\|\overline{s}_I\|_\infty\over\underline{M}^{1-{r\over 2}}}\|x_I^k\|^{1-{r\over 2}}_\infty\left\|X_{k,I}^{1-{r\over 2}}s^k_I\right\|_2\hfill (2bis)$ 

on the other hand. Since $\lim\limits_{k\uparrow\infty}x_I^k=0$, by $(2bis)$ we have necessarily $\left\|X_{k}^{1-{r\over 2}}s^k\right\|_\infty=\left\|X_{k,I}^{1-{r\over 2}}s^k_I\right\|_\infty={x_{i_{r\over 2}}^k}^{1-{r\over 2}}s_{i_{r\over 2}}^k$, for $k$ large enough. 
Now $\left\|X_{k,J}^{1-{r}}s_J^k\right\|_2\leq\underline{M}^{-{r\over 2}}\left\|X_{k,J}^{1-{r\over 2}}s_J^k\right\|_2$ 
and $\left\|X_{k,I}^{1-{r\over 2}}s^k_I\right\|_2=\left\|X_{k,I}^{{r\over 2}}X_{k,I}^{1-r}s^k_I\right\|_2\leq\|x_I^k\|_\infty^{r\over 2}\left\|X_{k,I}^{1-r}s^k_I\right\|_2$.
 Then using $(29bis)$ we get $\left\|X_{k,J}^{1-r}s_J^k\right\|_2\leq\displaystyle{L(A,c)\sqrt{n_I}\|\overline{s}_I\|_\infty\over\underline{M}}\|x_I^k\|_\infty\left\|X_{k,I}^{1-{r}}s^k_I\right\|_2.$
Hence using again the fact that $\lim\limits_{k\uparrow\infty}x_I^k=0$ we get $\left\|X_{k}^{1-{r}}s^k\right\|_\infty=\left\|X_{k,I}^{1-{r}}s^k_I\right\|_\infty$ for $k$ large enough.

\noindent Turn back now to (2). Then when $k$ is large enough we have

$\displaystyle{{L(A,c)}^2\over n_I^{1-{r\over 2}}}\|x^k_I\|_2^{2-r}\leq 2\|X_{k,I}^{1-{r\over 2}}s_I^k\|_2^2\leq2n_I\|X_{k,I}^{1-{r\over 2}}s_I^k\|_\infty^2$

$\hfill =2n_I\left({x^k_{i_{r\over 2}}}^{1-{r\over 2}}s_{i_{r\over 2}}^k\right)^2\leq 2n_I\left\|x^k_I\right\|_2^{2-r}\left|s^k_{i_{r\over 2}}\right|^2\hfill (3)$

\noindent and

$\displaystyle{{L(A,c)}^2\over n_I^{1-{r\over 2}}}\|x^k_I\|_2^{2-r}\leq 2n_I\left({x^k_{i_{r\over 2}}}^{1-{r\over 2}}s_{i_{r\over 2}}^k\right)^2\leq 2n_I\left({x^k_{i_{r\over 2}}}\right)^{2-r}\left\|s^k_I\right\|_\infty^2.\hfill (4)$

\noindent Using (3) and Lemma \ref{prelim1.1} it follows that

$\hfill\displaystyle\left|s^k_{i_r}\right|\geq\left|s^k_{i_{r\over 2}}\right|\geq{{L(A,c)}\over \sqrt{2}n_I^{1-{r\over 4}}}\hfill (5)$

\noindent and then $\left|s^k_{i_r}\right|=O(1)$. Now using  (4), (5) and the fact that $+\infty>g=\sup\limits_{k\in\mathbb{N}}\|s^k\|_\infty$, we get

$\hfill \left\|x_I^k\right\|_2\geq\displaystyle x^k_{i_{r\over 2}}\geq {L(A,c)^{2\over 2-r}\over 2^{1\over 2-r}g^{2\over 2-r}n_I^{1+{r\over 4-2r}}}\left\|x_I^k\right\|_2\hfill (6)$

\noindent and then $x^k_{i_{r\over 2}}=O\left(\left\|x_I^k\right\|_2\right)$ and then $i_{r\over 2}\in I$. Now $x_{i_r}^{k^{1-r}}O(1)=x_{i_r}^{k^{1-r}}\left|s^k_{i_r}\right|=\max\limits_i \left(x_i^{k^{1-r}}\left|s^k_i\right|\right)\geq x_{i_{r\over 2}}^{k^{1-r}}\left|s^k_{i_{r\over 2}}\right|=O\left(\left\|x_I^k\right\|_2\right)^{1-r}$. It follows that $\|X^{1-r}_{k,I}s^k_I\|_\infty=O(\|x^k_I\|^{1-r}_\infty)$ and $x_{i_r}^k=O\left(\left\|x_I^k\right\|_2\right)$, witch implies that $i_r\in I$.

Now using (5), (6) and the fact that $|s_{i_r}^k|\leq\|s^k\|_\infty\leq g$ we get
$\displaystyle x^k_{i_{r\over 2}}{\left|s^k_{i_{r\over 2}}\right|^2\over \left|s^k_{i_{r}}\right|}\geq {L(A,c)^{6-r\over 2-r}\over 2^{3-r\over 2-r}g^{4-r\over 2-r}n_I^{3-{r(1-r)\over 2(2-r)}}}\left\|x_I^k\right\|_2.$
But
$\langle c,x^k\rangle-{\overline x}=\langle {\overline s},x^k-{\overline x}\rangle=\langle {\overline s}_I,x^k_I\rangle\leq \|{\overline s}_I\|_2\|x^k_I\|_2\leq\sqrt{n_I}g\|x^k_I\|_2.$
It follows that
$\displaystyle x^k_{i_{r\over 2}}{\left|s^k_{i_{r\over 2}}\right|^2\over \left|s^k_{i_{r}}\right|}\geq {L(A,c)^{6-r\over 2-r}\over 2^{3-r\over 2-r}g^{6-2r\over 2-r}n_I^{{7\over 2}-{r(1-r)\over 2(2-r)}}}(\langle c,x^k\rangle-\overline{c})$
and then by (1), 
$\left\langle c,x^{k+1}\right\rangle -{\overline c}\leq{\overline L}(\langle c,x^k\rangle-{\overline c}).$

{\bf iii)} By { ii)} and Proposition \ref{convergencex},
$\displaystyle\|x^k_I\|_\infty\leq \left\|x^k-{\overline x}\right\|_2\leq H\left(\left\langle c,x^k\right\rangle -{\overline c}\right)
\leq H\overline{L}^{k-K}\left(\left\langle c,x^{K}\right\rangle-{\overline c}\right),\ \forall k\geq K.$
The result follows since $$\displaystyle0<\overline{L}^a=\left(1-{L(A,c)^{6-r\over 2-r}\over 2^{3-r\over 2-r}g^{6-r\over 2-r}n_I^{{7\over 2}-{r(1-r)\over 2(2-r)}}}\right)^a<1,\ \forall a>0.$$

{\bf iv)} We have $L(A,c)\|x^k_I\|_2\leq L(A,c)\|x^k-\overline{x}\|_2\leq\langle c,x^k\rangle-\overline{c}=\langle s^k,x^k-\overline{x}\rangle=\langle s^k_I,x^k_I\rangle+\langle s^k_J,x^k_J-\overline{x}_J\rangle=\langle \overline{s}_I,x^k_I\rangle\leq\|\overline{s}_I\|_2\|x^k_I\|_2$. The result then follows from  i).

{\bf v)} Let $({\overline y},{\overline s})$ be an accumulation point to $(y^k,s^k)$. We have $\langle c,x^k\rangle-{\overline c}=\langle A^t{\overline y}+{\overline s},x^k-{\overline x}\rangle=\langle {\overline s}_I,x^k_I\rangle$. Using Proposition \ref{convergencex} and the Cauchy-Schwartz inequality we get
$$\|{\overline s}\|_2\|x^k_I\|_2\geq\langle\overline{s},x^k\rangle=\langle c,x^k\rangle-{\overline c}\geq\displaystyle{1\over h}\|x^k-{\overline x}\|_2\geq{1\over h}\|x^k_I\|_2.$$
The result then follows.
\hfill\qed

Now we establish some technical results given by Saigal \cite{saigal} in the classical case.

\begin{proposition}\label{résultattechniques}
Let $(u^k)$ the sequence defined by
$u^k=\displaystyle {X_k^{r\over 2}PX_k^{1-{r\over 2}}c\over \left\langle c,x^k\right\rangle-{\overline c}}=\displaystyle {X_ks^k\over \left\langle c,x^k\right\rangle-{\overline c}}.$
Then we have:

{\bf i)} The sequence $(u^k)$ is bounded.

{\bf ii)} $\sum\limits_{k=0}^{\atop \infty}\left\|u_J^k\right\|^2<+\infty$.

{\bf iii)} $\sum\limits_{k=0}^{\atop \infty}|\delta^k|<+\infty$, where $\delta^k=\left\langle e_I,u^k_I\right\rangle -1$.

\end{proposition}

{\it Proof.}

{\bf i)} and {\bf ii)} By Proposition \ref{convergencex},  there is $h>0$ such that
$\left\|x^k-{\overline x}\right\|\leq h\left(\left\langle c,x^k\right\rangle -{\overline c}\right).$ Then $\|u_I\|_2=\displaystyle{\|X_{k,I}s^k_I\|_2\over \langle c,x^k\rangle-{\overline c}}\leq h{\|X_{k,I}s^k_I\|_2\over\|x^k-{\overline x}\|_2}\leq h{\|x_I^k\|_\infty\|s^k_I\|_2\over\|x^k_I\|_\infty}=h\|s^k_I\|_2$
Hence $u_I^k$ is bounded according to Proposition \ref{yksk}. According to Proposition \ref{convergencecx}
$\|u^k_J\|_2=\displaystyle{\|X_{k,J}s^k_J\|_2\over \langle c,x^k\rangle-{\overline c}}\leq h\|x^k_J\|_\infty{\|s^k_J\|_2\over\|x^k-{\overline x}\|_2}\leq h\overline{M}{\|s^k_J\|_2\over\|x^k_I\|_2}\leq h\hat{C}\overline{M}\|x_I^k\|^{1-r}.$
Then i) and ii) follows by using Proposition \ref{convergencecx}.


{\bf iii)} Set $SD=\left\{(y,s):\ A^ty+s=c,\ s_J=0\right\}$ the expected dual optimal solutions' set. Let $\left({\hat y}^k,{\hat s}^k\right)$
a solution to the problem $\min\left\{\left\|s^k-s\right\|_2:\ (y,s)\in SD\right\}.$ We have
$\left\langle c,x^k\right\rangle -{\overline c}=\left\langle c,x^k-{\overline x}\right\rangle=\left\langle {\hat s}^k+A^t{\hat y}^k,x^k-{\overline x}\right\rangle
=\left\langle {\hat s}_I^k,x_I^k\right\rangle.$ 
Then $\left|\left\langle e_I,u_I^k\right\rangle-1\right|
=\displaystyle{\left\langle x_I^k,s_I^k\right\rangle-\left\langle {\hat s}_I^k,x_I^k\right\rangle\over \left\langle c,x^k\right\rangle -{\overline c}}
\leq {\left\|s_I^k-{\hat s}_I^k\right\|_2\left\|x_I^k\right\|_2\over \left\langle c,x^k\right\rangle -{\overline c}}.$ By Proposition \ref{convergencex} we have
$\left\|x_I^k\right\|_2\leq \left\|x^k-{\overline x}\right\|_2\leq h\left(\left\langle c,x^k\right\rangle -{\overline c}\right).$ Hence
$\left|\left\langle e_I,u_I^k\right\rangle-1\right|_2\leq h\left\|s_I^k-{\hat s}_I^k\right\|_2.$
From Theorem 7 of \cite{saigal}, there is ${\hat M}$ such that $\left\|{\hat s}^k-s^k\right\|_2\leq {\hat M}\left\|s^k_J\right\|_2.$ Using Proposition \ref{convergencecx} we get $\left\|s_I^k-{\hat s}_I^k\right\|_2\leq {\hat M}\left\|s^k_J\right\|_2\leq {\hat C}{\hat M}\|x_I\|_\infty^{2-r}.$
The result then follows by using  iii) of Proposition \ref{convergencecx}.\hfill\qed

Now let us introduce the potential function $F$ defined as follow: $$F(x)=\ln\left(\displaystyle{\langle c,x\rangle-{\overline c}\over \tilde{\xi}_r(x)}\right)\mbox{, where }{\tilde \xi_r}:x\mapsto\left\{\begin{array}{cc}\left(\sum\limits_{i\in I}x_i^r\right)^{1\over r}&\mbox{if }x\geq 0,\cr
-\infty&\mbox{elsewhere.}\end{array}\right.$$  The following Proposition holds.

\begin{proposition}\label{potentialfunction1} There is $\Delta\in\mathbb{R}$ such that for every $k\geq0$, $F(x^k)\geq\Delta>-\infty.$
\end{proposition}

{\it Proof.}
By Theorem 3.1 and Proposition 4.3 of \cite{bacr}, ${\tilde\xi_r}(x^k){\tilde\xi_t}(e)\leq\langle x^k_I,e_I\rangle$, where $t$ is such that $\displaystyle{1\over t}+{1\over r}=1$. Hence ${\tilde\xi_r}(x^k)\leq\displaystyle{1\over n_I^{1\over t}}\sum\limits_{i\in I}x^k_i\leq{n^{1\over 2}\over n_I^{1\over t}}\left\|x_I^k\right\|_2=n^{{1\over r}-{1\over 2}} \left\|x_I^k\right\|_2=n^{{1\over r}-{1\over 2}}\left\|x^k-{\overline x}\right\|_2\leq n^{{1\over r}-{1\over 2}}h\left(\left\langle c,x^k\right\rangle-{\overline c}\right)$ and then $-\infty<\ln\left(\displaystyle{1\over n^{{1\over r}-{1\over 2}}h}\right)\leq F(x^k)$. The result then follows.\hfill\qed

\begin{proposition}\label{potentialfunction2}

Let $\displaystyle\beta\in\left(0,{2\over 3}\right)$. There exists $K\in\mathbb{N}$ such that $\forall k\geq K$
$$F(x^{k+1})-F(x^k)\leq-\theta^k\left(\Upsilon^k+(1-\Upsilon^k)\displaystyle{2-3\beta\over3(1-\beta)}\right)\|X_{k,I}^{r\over 2}w_I^k\|_2^2-(1+\Upsilon^k)\displaystyle{\theta^k\over\sum\limits_{j\in I}{x^k}^r}\delta^k-\theta^k\gamma^k$$
where $\Upsilon^k=\left\{\begin{array}{ll}1&\mbox{if }\max\limits_{i\in I}w_i^k\leq0,\cr
0&\mbox{if }\max\limits_{i\in I}w_i^k>0,\end{array}\right.$ $
\theta^k=\displaystyle{{\tilde t}^k\over 1-\displaystyle{1\over \sum\limits_{i\in I}{x_i^k}^r}{\tilde t}^k}$,
${\tilde t}^k=t^k\left(\left\langle c,x^k\right\rangle -{\overline c}\right)$, $t^k=\beta t^k_{\max}$,
$\gamma^k=\left\|X_{k,I}^{-{r\over 2}}u_J^k\right\|^2$ and
$\displaystyle w^k_I=X_{k,I}^{-r}u^k_I-{1\over \sum\limits_{i\in I}{x_i^k}^r}e$.
\end{proposition}

{\it Proof.}
Let us proof at first that given $\beta\in(0,1)$, there is $K\in\mathbb{N}$ such that $\forall k\geq K$, 
$F\left(x^{k+1}\right)-F\left(x^k\right)\leq\ln\left(1-\theta^k\| X_{k,I}^{r\over 2}w_I^k\|-\displaystyle2{\theta^k\over\sum\limits_{i\in I}{x^k_i}^r}\delta^k-\theta^k\gamma^k\right)-\displaystyle\sum\limits_{i\in I}\displaystyle{{x^k_i}^r\over\sum\limits_{i\in I}{x^k_i}^r}\ln\left(1-\theta^kw_i^k\right).$
We have $F(x^{k+1})-F(x^k)=\ln\left(\displaystyle{\langle c,x^{k+1}\rangle-{\overline c}\over\langle c,x^{k}\rangle-{\overline c}}\right)-{1\over r}\ln\left(\displaystyle{\sum x_i^{{k+1}^r}\over\sum x_i^{{k}^r}}\right)$
,
$\langle c,x^{k+1}\rangle-{\overline c}=\langle c,x^{k}\rangle-{\overline c}-t^k\langle X_k^{-r}X_kS^k,X_ks^k\rangle$,
 $\displaystyle u^k={X_ks^k\over \langle c,x^{k}\rangle-{\overline c}}$ and ${\tilde t^k}=(\langle c,x^{k}\rangle-{\overline c})t^k$. Then
$\displaystyle{\langle c,x^{k+1}\rangle-{\overline c}\over \langle c,x^{k}\rangle-{\overline c}}=1-{\tilde t^k}\langle X^{-r}_{k,I}u^k,u^k\rangle=1-{\tilde t^k}\langle X^{-r}_{k,I}u^k_I,u^k_I\rangle -\tilde{t^k}\gamma^k.$
Now $u_I^k=X_{k,I}^rw_I^k+{1\over \sum\limits_{i\in I}x_i^{k^r}}X_{k,I}^re_I$, $X_{k,I}^{-r}u_I^k=w_I^k+{1\over \sum\limits_{i\in I}x_i^{k^r}}e_I$ and $\delta^k=\langle u^k_I,e_I\rangle-1=\langle X^r_{k,I}w_I^k,e_I\rangle$. Then $\displaystyle \langle X^{-r}_{k,I}u^k_I,u^k_I\rangle=\langle X^{r}_{k,I}w^k_I,w^k_I\rangle+{2\over \sum\limits_{i\in I}x_i^{k^r}}\delta^k+{1\over \sum\limits_{i\in I}x_i^{k^r}}$. It follows that

$\begin{array}{ll}\displaystyle{\langle c,x^{k+1}\rangle-{\overline c}\over \langle c,x^{k}\rangle-{\overline c}}&=\displaystyle 1-{\tilde t^k}\langle X_{k,I}^rw_I^k,w_I^k\rangle-{2{\tilde t^k}\over \sum\limits_{i\in I}x_i^{k^r}}\delta^k-{{\tilde t^k}\over \sum\limits_{i\in I}x_i^{k^r}}-{\tilde t^k}\gamma^k\cr
&=\displaystyle\left(1-{{\tilde t^k}\over \sum\limits_{i\in I}x_i^{k^r}}\right)\left[1-\theta^k\langle X_{k,I}^rw_I^k,w_I^k\rangle-{2\theta^k\over \sum\limits_{i\in I}x_i^{k^r}}\delta^k-\theta^k\gamma^k\right].\end{array}$

Let us show now, for $\beta\in (0,1)$,
$\displaystyle\theta^k={{\tilde t^k}\over 1-{{\tilde t^k}\over \sum\limits_{i\in I}x_i^{k^r}}}={\sum\limits_{i\in I}x_i^{k^r}{\tilde t^k}\over \sum\limits_{i\in I}x_i^{k^r}-{\tilde t^k}}>0$, for $k$ large enough.
We have $\langle c,x^k\rangle-{\overline c}=\langle s^k,x^k-\overline{x}\rangle=\langle s^k_I,x^k_I\rangle+\langle s^k_J,x^k_J-\overline{x}_J\rangle$.

$\begin{array}{ll}\mbox{Then }\sum\limits_{i\in I}x_i^{k^r}-{\tilde t^k}&=\sum\limits_{i\in I}x_i^{k^r}-\displaystyle\beta{\langle c,x^k\rangle-{\overline c}\rangle\over \max\limits_{i\in I}(X_k^{1-r}s^k)_i}\cr&=\displaystyle{\sum\limits_{i\in I}x_i^{k^r}\max\limits_{i\in I}(X_k^{1-r}s^k)_i-\beta(\langle s^k_I,x^k_I\rangle+\langle s^k_J,x^k_J-\overline{x}_J\rangle)\over \max\limits_{i\in I}(X_k^{1-r}s^k)_i}\cr
&\geq\displaystyle{\sum\limits_{i\in I}x_i^{k^r}(x_i^{k^{1-r}}s^k_i)-\beta\langle s^k_I,x^k_I\rangle-\beta\langle s^k_J,x^k_J-\overline{x}_J\rangle)\over \max\limits_{i\in I}(X_k^{1-r}s^k)_i}\cr
&=\displaystyle{\langle x_I^k,s_I^k\rangle\over\max\limits_{i\in I}(X_k^{1-r}s^k)_i}\left(1-\beta\displaystyle{\langle s_J^k,x_J^k-\overline{x}_J\rangle\over\langle x_I^k,s_I^k\rangle}-\beta\right)\cr
&\geq \displaystyle{\langle x_I^k,s_I^k\rangle\over\max\limits_{i\in I}(X_k^{1-r}s^k)_i}\left(1-\beta\displaystyle{ \|s_J^k\|_2\|x_J^k-\overline{x}_J\|_2\over\langle x_I^k,s_I^k\rangle}-\beta\right)\cr
\end{array}
$

Using the fact that $\lim\limits_{k\uparrow\infty}x^k-\overline{x}=0$, from Proposition \ref{convergencecx} we get
$\displaystyle{ \|s_J^k\|_2\|x_J^k-\overline{x}_J\|_2\over\langle x_I^k,s_I^k\rangle}=o(\|x_I^k\|^{1-r})$. Therefore $\theta^k>0$.
Now $\displaystyle{1\over r}\ln\left(\displaystyle{\sum\limits_{i\in I}x_i^{k^{r+1}}\over \sum\limits_{i\in I}x_i^{k^r}}\right)={1\over r}\ln\left(\sum\limits_{i\in I}\displaystyle{x_i^{k^r}\over\sum\limits_{j\in I}x_j^{k^r}}(1-{\tilde t^k}x_i^{k^{-r}}u_i^k)^r\right).$ Since the function $t\mapsto \ln t$, $t>0$, is concave, one has

$\begin{array}{ll}\displaystyle{1\over r}\ln\left(\displaystyle{\sum\limits_{i\in I}x_i^{k^{r+1}}\over \sum\limits_{i\in I}x_i^{k^r}}\right)^r&\geq\sum\limits_{i\in I}\displaystyle{x_i^{k^r}\over\sum\limits_{j\in I}x_j^{k^r}}\ln\left(1-{\tilde t^k}x_i^{k^{-r}}u_i^k\right)\cr
&=\sum\limits_{i\in I}\displaystyle{x_i^{k^r}\over\sum\limits_{j\in I}x_j^{k^r}}\ln\left(1-{{\tilde t^k}\over \sum\limits_{j\in I}x_j^{k^r}}-{\tilde t^k}w_i^k\right)\cr
&=\displaystyle\ln\left(1-{{\tilde t^k}\over \sum\limits_{j\in I}x_j^{k^r}}\right)+\sum\limits_{i\in I}\displaystyle{x_i^{k^r}\over\sum\limits_{j\in I}x_j^{k^r}}\ln\left(1-\theta^kw_i^k\right).
\end{array}$

Hence
$F(x^{k+1})-F(x^k)\leq\ln\left[1-\theta^k\langle X_{k,I}^rw_I^k,w_I^k\rangle-{2\theta^k\over \sum\limits_{i\in I}x_i^{k^r}}\delta^k-\theta^k\gamma^k\right]-\sum\limits_{i\in I}\displaystyle{x_i^{k^r}\over\sum\limits_{j\in I}x_j^{k^r}}\ln\left(1-\theta^kw_i^k\right).$
Assume now that $\max\limits_{i\in I}w_i^k>0$. Then using Lemma 8 of \cite{saigal} or its proof we can easily see
$\sum\limits_{i\in I}\displaystyle{x_i^{k^r}\over\sum\limits_{j\in I}x_j^{k^r}}\ln\left(1-\theta^kw_i^k\right)\geq\displaystyle-{\theta^k\over\sum\limits_{j\in I}x_j^{k^r}}\delta^k-\displaystyle{{\theta^k}^2\over 2\sum\limits_{j\in I}x_j^{k^r}}{\|X_{k,I}^{r\over 2}w_I^k\|^2\over 1-\theta^k\max\limits_{i\in I}w_i^k}$.
Using in addition the fact that
$\ln(1-a)\leq-a,\ \forall a<1$
we get
$F(x^{k+1})-F(x^k)\leq-\theta^k\left(1-\displaystyle{{\theta^k}\over 2\sum\limits_{j\in I}x_j^{k^r}}{1\over 1-\theta^k\max\limits_{i\in I}w_i^k}\right)\|X_{k,I}^{r\over 2}w_I^k\|_2^2-\displaystyle{{\theta^k}\over \sum\limits_{j\in I}x_j^{k^r}}\delta^k-\theta^k\gamma^k.$
Now 

$\begin{array}{ll}1-\displaystyle{{\theta^k}\over 2\sum\limits_{j\in I}x_j^{k^r}}{1\over 1-\theta^k\max\limits_{i\in I}w_i^k}&=1-\displaystyle{{\theta^k}\over 2\sum\limits_{j\in I}x_j^{k^r}}\displaystyle{1\over 1+\displaystyle{{\theta^k}\over \sum\limits_{j\in I}x_j^{k^r}}-\theta^k\max\limits_{i\in I}(X_k^{-r}u^k)_i}\cr
&=1-\displaystyle{1\over 2\sum\limits_{j\in I}x_i^{k^r}}\displaystyle{\theta^k\over 1+\displaystyle{{\theta^k}\over \sum\limits_{j\in I}x_j^{k^r}}}\displaystyle{1\over 1-\displaystyle{\theta^k\over 1+\displaystyle{{\theta^k}\over \sum\limits_{j\in I}x_j^{k^r}}}\max\limits_{i\in I}(X_k^{-r}u^k)_i}.
\end{array}
$

We have on the one hand $\theta^k=\displaystyle{{\tilde t^k}\over 1-\displaystyle{{\tilde t^k}\over \sum\limits_{j\in I}x_j^{k^r}}}$ and then ${\tilde t^k}=\displaystyle{\theta^k\over 1+\displaystyle{\theta^k\over \sum\limits_{j\in I}x_j^{k^r}}}$. On the other hand, ${\tilde t^k}=\beta\displaystyle{\langle c,x^k\rangle-{\overline c}\over \max\limits_{i\in I}(X_k^{1-r}s^k)_i}=\displaystyle{\beta\over \max\limits_{i\in I}(X_k^{-r}u^k)_i}$.
It follows that

$\begin{array}{ll}1-\displaystyle{{\theta^k}\over 2\sum\limits_{j\in I}x_j^{k^r}}{1\over 1-\theta^k\max\limits_{i\in I}w_i^k}&=1-\displaystyle{{\tilde t^k}\over 2(1-\beta)\sum\limits_{i\in I}x_i^{k^r}}\cr
&=1-\displaystyle{\beta\over 2(1-\beta)}\displaystyle{\langle c,x^k\rangle-{\overline c}\over\sum\limits_{i\in I}\left[x_i^{k^r}\max\limits_{j\in I}(X_k^{1-r}s^k)_i\right]}.
\end{array}
$

Using the fact that $\sum\limits_{i\in I}\left[x_i^{k^r}\max\limits_{j\in I}(X_k^{1-r}s^k)_i\right]\geq\sum\limits_{i\in I}{x_i^k}^r{x_i^k}^{1-r}s_i^k=\sum\limits_{i\in I}x_i^ks_i^k$ and the fact that $\langle c,x^k\rangle-{\overline c}=\langle s^k,x^k-\overline{x}\rangle=\langle s^k_I,x^k_I\rangle+\langle s_J^k,x_J^k-\overline{x}_J\rangle$ we get
$1-\displaystyle{\beta\over 2(1-\beta)}\displaystyle{\langle c,x^k\rangle-{\overline c}\over\sum\limits_{i\in I}\left[x_i^{k^r}\max\limits_{j\in I}(X_k^{1-r}s^k)_i\right]}\geq 1-\displaystyle{\beta\over 2(1-\beta)}\left(1+\displaystyle{\langle x_J^k,s_J^k\rangle\over\langle s_I^k,x_I^k\rangle}\right)
=\displaystyle{2-3\beta\over 2(1-\beta)}-{\beta\over2(1-\beta)}{\langle s_J^k,x^k_J-\overline{x}_J\rangle\over\langle x_I^k,s^k_I\rangle}
\geq\displaystyle{2-3\beta\over 2(1-\beta)}-{\beta\over2(1-\beta)}{\|s_J^k\|_2\|x^k_J-\overline{x}_J\|_2\over\langle x_I^k,s^k_I\rangle}
=\displaystyle{2-3\beta\over 2(1-\beta)}+o(\|x_I^k\|^{1-r}_2)
\geq\displaystyle{2-3\beta\over 3(1-\beta)},$ for k being large enough.
Hence
$F(x^{k+1})-F(x^k)\leq-\theta^k\displaystyle{2-3\beta\over3(1-\beta)}\|X_{k,I}^{r\over 2}w_I^k\|_2^2-\displaystyle{\theta^k\over\sum\limits_{j\in I}{x^k}^r}\delta^k-\theta^k\gamma^k.$
Consider now the case $\max\limits_{i\in I}w_i^k\leq0$. Then $\displaystyle\sum\limits_{i\in I}\displaystyle{{x^k_i}^r\over\sum\limits_{i\in I}{x^k_i}^r}\ln\left(1-\theta^kw_i^k\right)\geq0$ and the result follows by using the fact that $\ln(1-a)\leq-a,\ \forall a<1$.\hfill\qed
\begin{lemma}\label{lemme2}
We have
$\displaystyle{\theta^k\over\sum\limits_{j\in I}{x_i^k}^r}=\displaystyle{\displaystyle{{\tilde t^k}\over\sum\limits_{j\in I}{x_i^k}^r}\over1-\displaystyle{{\tilde t^k}\over\sum\limits_{j\in I}{x_i^k}^r}}=O(1),$
For $k$ large enough.
\end{lemma}

{\it Proof.} We have $n_Ig\|x^k_I\|_\infty\geq\sum\limits_{i\in I}{x_i^k}^r\max\limits_{i\in I}(X_k^{1-r}s^k)_i\geq\sum\limits_{i\in I}{x_i^k}^r(X_k^{1-r}s^k)_i=\langle x_I^k,s^k_I\rangle$. By Proposition \ref{convergencecx} $\langle x_I^k,s^k_I\rangle=O(\|x^k_I\|_\infty)$ and $\langle c,x^k\rangle-\overline{c}=O(\|x^k_I\|_\infty)$. It follows that $\sum\limits_{i\in I}{x_i^k}^r\max\limits_{i\in I}(X_k^{1-r}s^k)_i=O(\|x^k_I\|_\infty)$ and then  $\displaystyle{{\tilde t^k}\over \sum\limits_{i\in I}{x_i^k}^r}=\beta\displaystyle{\langle c,x^k\rangle -{\overline c}\over \sum\limits_{i\in I}{x_i^k}^r\max\limits_{i\in I}(X_k^{1-r}s^k)_i}=O(1)$. Now $\langle c,x^k\rangle-{\overline c}=\langle s^k,x^k-\overline{x}\rangle=\langle s^k_I,x^k_I\rangle+\langle s^k_J,x^k_J-\overline{x}_J\rangle\leq\langle s^k_I,x^k_I\rangle+\| s^k_J\|_2\|x^k_J-\overline{x}_J\|_2$. Then using again Proposition \ref{convergencecx} we get $\displaystyle{{\tilde t^k}\over \sum\limits_{i\in I}{x_i^k}^r}\leq\beta+\beta\displaystyle{\| s^k_J\|_2\|x^k_J-\overline{x}_J\|_2\over\langle x^k_I,s^k_I\rangle}=\beta+o\left(\|x_I^k\|^{1-r}\right).$ But $\beta\in(0,1)$. The result then follows.\hfill\qed

Now as is mentioned in Theorem 1 of \cite{ba}, given $t$ in $(-\infty,0)\cup(0,1)$, we say the $\xi_t$-dual-analytic center the unique optimal solution to the problem
$$\max\limits_{y,s}\{\xi_{t,I}(s_I):\ A^ty+s=c,\ s_J=0,\ s\ge0\}$$
where
$$\xi_{t,I}(x)=\left\{\begin{array}{ll}\left(\sum\limits_{i\in I}x_i^t\right)^{1\over t}&\mbox{if }x_I\geq0\cr
-\infty&\mbox{elsewhere,}\end{array}\right.$$ if $t\in(0,1)$ and $$\xi_{t,I}(s_I)=\left\{\begin{array}{ll}\left(\sum\limits_{i\in I}{s_i}^t\right)^{1\over t}&\mbox{if }s_I\in(0,+\infty)^{n_I}\cr 0&\mbox{if }s_I\in\partial([0,+\infty)^{n_I})\cr-\infty&\mbox{elsewhere},\end{array}\right.$$ if $t\in(-\infty,0)$.
The unicity of the solution is ensured by the strict quasi-concavity of $\xi_t$ and Lemma 1 of \cite{ba}. The KKT optimality conditions are then expressed as follow. There exist $(y,s)\in\mathbb{R}^m\times\mathbb{R}^n$ and $v\in\mathbb{R}^n$ satisfying the following conditions, $\nabla\xi_{t,I}(s_I)=v_I$, $Av=0$, $A^ty+s=c$, $s_J=0$, $s\geq0$.


{\it Proof of Theorem \ref{convergenceys}.}

According to Proposition \ref{potentialfunction1}, $\sum\limits_{k\geq0}(F(x^{k+1})-F(x^k))$ converges and according to Proposition \ref{résultattechniques}, $\sum\limits_{k\geq0}\delta^k$ and $\sum\limits_{k\geq0}\gamma^k$ converge too. Then using Proposition \ref{potentialfunction2}, $\sum\limits_{k\geq0}\|X_{k,I}^{r\over 2}w_I^k\|_2^2$ converges. Hence

$\hfill\displaystyle\lim\limits_{k\to +\infty}{X_{k,I}^{r\over2}\over\sum\limits_{i\in I} {x_i^k}^r}\left({\sum\limits_{i\in I} {x_i^k}^r\over \langle c,x^k\rangle-{\overline c}}X_{k,I}^{1-r}s^k_I-e\right)=0\hfill (1)$

For all $k\in \mathbb{N}$ we set $I(k)=\{i:\ \displaystyle{x_i^k\geq\|x^k_I\|_\infty^2}\}$. We shall prove that for some $K$ chosen large enough, $I(k)\subset I(k+1)$, $\forall k\geq K$. For $k\in\mathbb{N}$ and $i\in\{1,\cdots,n\}$ we set $\epsilon^k_i=\displaystyle{\sum\limits_{j\in I} {x_j^{k}}^r\over \langle c,x^{k}\rangle-{\overline c}}{x_i^{k}}^{1-r}s_i^{k}-1$. Let $\epsilon>0$ be small enough. Then by (1) there exists $K\in\mathbb{N}$ large enough such that $\forall k\geq K$, $|\epsilon_i^k|\leq\epsilon$, $\forall i\in I(k)$. Let $k\geq K$ and $i\in I(k)$. Since $K$ is assumed to be large we have necessarily from (1) $s_i^k>0$. Using in addition the fact that $\|s^k\|_\infty\leq g$ and Proposition \ref{convergencex} we get
$gn_I\|x^k_I\|^r_\infty{x_i^k}^{1-r}s_i^k\geq\sum\limits_{j\in I} {x_j^{k}}^r{x_i^{k}}^{1-r}s_i^{k}=(\langle c,x^{k}\rangle-{\overline c})(1+\epsilon_i^k)
\geq L(A,c)\|x^k-{\overline x}\|_2(1+\epsilon^k_i)\geq L(A,c)\|x^k_I\|_\infty(1-\epsilon).$
Hence $\|x^k_I\|_\infty\geq x_i^k\geq\displaystyle\left({L(A,c)(1-\epsilon)\over n_Ig}\right)^{1\over 1-r}\|x^k_I\|_\infty$
and then $x_i^k=O(\|x_I^k\|_\infty)$. Since in addition $x_i^{k+1}=\left(1-\beta\displaystyle{{x_i^k}^{1-r}s_i^k\over\max\limits_{j\in I}(X_k^{1-r}s^k)_j}\right)x_i^k$ and $0<\displaystyle{{x_i^k}^{1-r}s_i^k\over\max\limits_{j\in I}(X_k^{1-r}s^k)_j}\leq1$ we have $(1-\beta)x_i^k\leq x_i^{k+1}\leq (1+\beta)x_i^k$ and then $x_i^{k+1}=O(x_i^k)$. Now $\|x_I^{k+1}\|_2=\left\|x_I^k-\beta\displaystyle{{X_{k,I}}^{2-r}s_I^k\over\max\limits_{j\in I}(X_k^{1-r}s^k)_j}\right\|_2\leq\|x_I^k\|_2+\|x_I^k\|_2\displaystyle{\|{X_{k,I}}^{1-r}s_I^k\|_2\over\max\limits_{j\in I}(X_k^{1-r}s^k)_j}$ and $\sum\limits_{j\in I}{x_j^k}^r\|X_k^{1-r}s^k\|_\infty\geq\sum\limits_{j\in I}{x_j^k}^r\max\limits_{j\in I}{x_j^k}^{1-r}s_j^k\geq\sum\limits_{j\in I}{x_j^k}^r{x_j^k}^{1-r}s_j^k=\langle s_I,x_I\rangle.$ It follows that $\|X_k^{1-r}s^k\|_\infty\geq\max\limits_{j\in I}{x_j^k}^{1-r}s_j^k\geq\displaystyle{\langle s_I,x_I\rangle\over\sum\limits_{j\in I}{x_j^k}^r}.$
But (Proposition \ref{convergencecx}) $\|X_k^{1-r}s^k\|_\infty=O(\|x_I^k\|_\infty^{1-r})$, $\langle  s^k_I,x_I^k\rangle=O(\|x_I^k\|_\infty)$ and $\sum\limits_{j\in I}{x_j^k}^r=O(\|x_I^k\|_\infty^r)$. Then $\max\limits_{j\in I}{x_j^k}^{1-r}s_j^k=O(\|x_I^k\|_\infty^{1-r})$ and then $\displaystyle{\|{X_{k,I}}^{1-r}s_I^k\|_2\over\max\limits_{j\in I}(X_k^{1-r}s^k)_j}=O(1).$ Hence there is $\varrho>0$ such that $\|x_I^{k+1}\|_2\leq\varrho\|x_I^k\|_2$. So $x_i^{k+1}=O(x_i^k)=O(\|x_I^k\|_2)\geq O(\|x_I^{k+1}\|_2)\geq\|x_I^{k+1}\|_2^2$ and then $i\in I(k+1)$.
Set now ${\hat I}=\cup{\atop{k\in\mathbb{N}}}I(k)$ and let us prove that in fact ${\hat I}=I$. Assume for contradiction that there is $i\in I-{\hat I}$. Then $\forall k\geq K$, $\displaystyle{x_i^{k+1}\over x_i^k}=1-\beta{{x_i^{k}}^{1-r}s_i^{k}\over\|X_{I,k}^{1-r}s^{k}_I\|_\infty}=1-\beta{{x_i^{k}}^{1-r}\over\|x^k_I\|^{2(1-r)}_\infty}{\|x_I^k\|_\infty^{1-r}\over\|X_{I,k}^{1-r}s^{k}_I\|_\infty}s_i^{k}\|x_I^k\|_\infty^{1-r}\geq1-\beta{\|x_I^k\|_\infty^{1-r}\over\|X_{I,k}^{1-r}s^{k}_I\|_\infty}s_i^{k}\|x_I^k\|_\infty^{1-r}$. We know by Proposition \ref{convergencecx} that $\|X_{I,k}^{1-r}s^{k}_I\|_\infty=O(\|x_I^k\|_\infty^{1-r})$. Using in addition the fact that $s^k$ is bounded it follows that $\displaystyle{x_i^{k+1}\over x_i^k}\geq 1+O(\|x_I^k\|_\infty^{1-r})$. Hence for all $K'\geq K$, $\displaystyle{x_i^{K'}\over x_i^K}\geq \prod\limits_{K\leq k\leq K'}(1+O(\|x_I^k\|_\infty^{1-r}))=1+O\left(\sum\limits_{k=K}^{K'}\|x_I^k\|_\infty^{1-r}\right)$. We know by Proposition $\ref{convergencecx}$ that $\sum\limits_{k\in\mathbb{N}}\|x_I^k\|_\infty^{1-r}$ is a converging serie. It follows that chosing $K$ large enough, $O\left(\sum\limits_{k=K}^{+\infty}\|x_I^k\|_\infty^{1-r}\right)\geq\displaystyle-{1\over 2}$. Then 
$0=\lim\limits_{K'\to+\infty}\displaystyle{x_i^{K'}\over x_i^K}\geq{1\over 2}$, which is absurde. Hence

$\hfill\lim\limits_{k\to+\infty}\epsilon^k=\lim\limits_{k\to+\infty}\displaystyle{\sum\limits_{i\in I} {x_i^k}^r\over \langle c,x^k\rangle-{\overline c}}X_{k,I}^{1-r}s^k_I-e_I=0\hfill (3)$
 
and then there is necessarily $\tau>0$ such that $\tau e_I<s_I^k$ for $k$ being large enough. Let now $({\tilde y},{\tilde s})$ be an accumulation point to $(y^k,s^k)$. Then we have ${\overline X}{\tilde s}=0,\ A{\overline x}=b,\ A{\tilde y}+{\tilde s}=c,\ {\overline x}\geq0$ and ${\tilde s}\geq 0$. The KKT optimality conditions of $(LP)$ are then satisfied and then ${\overline x}$ is an $(LP)$ optimal solution and $({\tilde y},{\tilde s})$ is a dual optimal solution. Moreover since ${\overline x}_I=0,\ {\tilde s}_I>0,\ x_J>0$ and ${\tilde s}_J=0$ the strict complementary slackness condition holds.
Let now $k$ be large enough. Then it is easy to see with the help of Proposition \ref{convergencecx} that $(\langle c,x^k\rangle-{\overline c})\displaystyle{X_{k,I}^{r-1}\over\sum\limits_{i\in I} {x_i^k}^r}=O(1)$. It follows that 
$s_I^k=(\langle c,x^k\rangle-{\overline c})\displaystyle{X_{k,I}^{r-1}\over\sum\limits_{i\in I} {x_i^k}^r}(e_I+\epsilon^k)=(\langle c,x^k\rangle-{\overline c})\displaystyle{X_{k,I}^{r-1}e_I\over\sum\limits_{i\in I} {x_i^k}^r}+\hat{O}(\epsilon)
=\displaystyle{(\langle c,x^k\rangle-{\overline c})\over\xi_{r,I}(x_I^k)}\nabla\xi_{r,I}(x_I^k)+\hat{O}(\epsilon),$
where $\hat{O}({\epsilon})$ represents every function of $\epsilon$ satisfying $\lim\limits_{\epsilon\downarrow0}\hat{O}({\epsilon})=0$. According to Theorem 3.1 of \cite{bacr}, we have $\xi_{t,I}(\nabla\xi_{r,I}(x_I^k))=1$ and then, by continuity of $\xi_{t,I}$ on $(0,+\infty)^{n_I}$, we get $\xi_{t,I}(s^k_I)=\displaystyle{(\langle c,x^k\rangle-{\overline c})\over\xi_{r,I}(x_I^k)}+\hat{O}(\epsilon)=O(1)$. Hence $\displaystyle{s_I^k\over\xi_{t,I}(s^k_I)}=\nabla\xi_{r,I}(x_I^k)+\hat{O}(\epsilon).$ Now it is easy to see from Theorem 3.1 of \cite{bacr} that $\nabla\xi_{r,I}(\cdot)$ is positively homogeneous of degree 0. It follows that $\displaystyle{s_I^k\over\xi_{t,I}(s^k_I)}=\nabla\xi_{r,I}\left({x_I^k\over\|x_I^k\|_\infty}\right)+\hat{O}(\epsilon).$ By (3) there is $\tau'>0$ such that $\tau'e_I\leq\displaystyle{x_I^k\over\|x_I^k\|_\infty}\leq e_I$ for $k$ being large enough. Then using in addition Lemma \ref{prelim2.1}, $z^k=\displaystyle{x^k-\overline{x}\over\|x_I^k\|_\infty}$ is bounded. Let then $\overline{z}$ be a limit of a convergent subsequence of $(z^k)$. Then
$\displaystyle{\overline{s}_I\over\xi_{t,I}(\overline{s}_I)}=\nabla\xi_{r,I}(\overline{z}_I).$ Using again Theorem 3.1 of \cite{bacr} one has
$\displaystyle{\overline{z}_I\over\xi_{r,I}(\overline{z}_I)}=\nabla\xi_{t,I}(\overline{s}_I).$ But $Az^k=0$. It follows that $A\left(\displaystyle{\overline{z}\over\xi_{r,I}(\overline{z}_I)}\right)=0$. Hence $\left(\overline{y},\overline{s},\displaystyle{\overline{z}\over\xi_{r,I}(\overline{z}_I)}\right)$ satisfies the KKT optimality conditions for the problem
$$\max\limits_{y,s}\{\xi_{t,I}(s_I):\ A^ty+s=c,\ s_J=0,\ s\ge0\}$$
The result then follows\hfill\qed



Turn back now to the case where ${\cal I}\not=\emptyset$. Then using (\ref{kktldb}) of Section \ref{algo} and adapting results of this section, the expected dual approximate optimal solution vectors $y,\ s$ and $w$, associated to a current point $x$ are

\begin{eqnarray}\label{dualy}y=(AH^{-{1\over2}}A^t)^{-1}H^{-{1\over2}}c
\end{eqnarray}
\begin{eqnarray}\label{dualw}
 w_{\cal I}=-U_{\cal I}^{-1}X_{{\cal I}}(c-A^ty)_{\cal I}=-U_{\cal I}^{-1}X_{{\cal I}}(H^{{1\over2}}PH^{-{1\over2}}c)_{\cal I}\mbox{ and }w_{\cal\overline{I}}=0
\end{eqnarray}
\begin{eqnarray}\label{duals}
 s=c-A^ty+w=H^{{1\over2}}PH^{-{1\over2}}c+w
\end{eqnarray}
where $U=diag(u)$. Hence the expected relative duality gap is

\begin{eqnarray}\label{rgap}
Rgap=\displaystyle{\langle c,x\rangle-\langle b,y\rangle+\langle u_{\cal I},w_{\cal I}\rangle\over|\langle c,x\rangle|+1}
\end{eqnarray}
%
%
%

\section{Numerical results}\label{numerical}

The porpose of the following tests is to compare the algorithm performance according to different values of $r$ between 0 and 1. We have opted to consider the following values $r=0$ (the classical case), $r=0.1$, $r=0.2$, $r=0.3$, $r=0.4$, $r=0.5$, $r=0.6$ and $r=0.7$. We solved a large set of testing problems, taken from the familiar Netlib test set (GAY \cite{netlib}). For values of $r$ exceeding 0.8 the algorithm showed lesser efficiency on most problems. 
Results obtained are shown in Table 1. Each row in the table contains the name of the problem and the number of iterations for the different values of $r$. The parameter of the stopping rule is $\epsilon=10^{-10}$. To read the mps-files, the specs-files and perform the symbolic Cholesky factorization we use rdmps1, rdmps2, rdspec, prepro and all dependencies written by Gondzio\cite{hopdm}. We use no presolving. Our numerical experiments were performed on a laptop {\it hp ZBook} (Processor: Intel core i7-4810 MQ, CPU $2.804\times 8$HZ - Operating system: Ubuntu linux). The code is written in GNU Fortran 95.

\setlength{\textwidth}{420pt}
\begin{longtable}{|l|l|l|l|l|l|l|l|l|}\hline
Problem   & r=0 &r=0.1 & r=0.2& r=0.3& r=0.4& r=0.5&r=0.6 & r=0.7\\\hline\hline
25fv47    & 58  & 59   &  60  & 64   &  66  & 81   & 90   & $\star\star$ \\
80bau3b   & 80 & 79  &  77 & 76   &  78  & 75   & 91   & $\star\star$  \\
adlittle  & 34  & 34   &  33  & 33   &  34  & 37   & 40   & 57   \\
afiro     & 25  & 23   &  23  & 21   &  21  & 20   & 21   & 22  \\
agg       & 45  & 49   &  43  & 62   &  66  & 85   & 129  & $\star\star$  \\
agg2      & 41  & 41   &  42  & 46   &  48  & 54   & 69   & $\star\star$  \\
agg3      & 39  & 40   &  43  & 43   &  48  & 56   & 69   & $\star\star$  \\
bandm     & 46  & 52   &  53  & 60   &  72  & 96   & 143   & $\star\star$  \\
beaconfd  & 35  & 34   &  33  & 32   &  30  & 30   & 30   & 32 \\
blend     & 40  & 41   &  40  & 43   &  44  & 37   & 44   & 49 \\
bnl1      & 66  & 67   &  60  & 65   &  75  & 92   &$\star\star$& $\star\star$ \\
bnl2      & 75  & 87   &  91  & 99   &  120  & 148  & $\star\star$   & $\star\star$ \\
boeing1   & 62  & 71   &  70  & 71   &  72  & 103  & 149  & $\star\star$ \\
boeing2   & 40  & 57   &  45  & 50   &  59  & 70   & 92  & 131 \\
bore3d    & 67  & 72   &  80  & 86   &  95  & 121  & 147  & $\star\star$ \\
brandy    & 36  & 37   &  38  & 39   &  44  & 53   & 66   & 92 \\
capri     & $\star\star$  & $\star\star$   &  $\star\star$  & $\star\star$   &  49  & $\star\star$   &$\star\star$   & $\star\star$  \\
cycle     &101&$\star\star$&90 & 107  &  118 & 136  &165&$\star\star$ \\
czprob    & 80  & 64   &  61  & 57   &  59  &  61 & 77   & $\star\star$ \\
d2q06c    & 71 &69&  69  & 72   &  75 &103 &134&$\star\star$  \\
d6cube    & 117 &87&  85  & 72   & 68 &70 &74&82  \\
degen2    & 44  & 43  &  43  & 37   &  42  & 44  & 49   & 62  \\
degen3    & 53  & 51   &  50  & 47   &  51  & 49   & 64   & 91  \\
dfl001    & $\star\star$&$\star\star$&  144 & 166 &  174& $\star\star$   & $\star\star$   & $\star\star$\\
e226      & 53  & 53   &  52  & 51   &  69  & 74   & 88   & $\star\star$   \\
etamacro  & 49  & 50   &  55  & 65   &  79  & 100   &$\star\star$ & $\star\star$   \\
fffff800  &61&56&55& 58&68&62& 140 &212 \\
finnis    &78 &82 &80&85&$\star\star$&$\star\star$& $\star\star$& $\star\star$ \\
fit1d     &41  & 40   & 66  &65  &  52  & 76   &$\star\star$   &$\star\star$\\
fit1p     &31& 32 &32  & 41   & 48  &50   & 73 & $\star\star$   \\
fit2d     &55 &45 &  43  &42   &  44  & 58   & $\star\star$   & $\star\star$   \\
fit2p     &44  &47  &48  &51   &45 &56  &73  &110   \\
forplan   & 48  & 46   & 46  &50   & 52  & 55   & 64   & $\star\star$   \\
ganges    & 27  &26   &34  &29   &32  & 36   &44   & 56   \\
gfrd-pnc  &29  &30   &37  & 49   & 65  & 89  & 134  & $\star\star$   \\
greenbea  &$\star\star$& $112^*$ &105  &  $99^*$ & $117^*$  & $132^*$ &$181^*$&$\star\star$ \\
greenbeb  & 85  & 81   &  84  & 90   &  98  & 112& $\star\star$&$\star\star$\\
grow7     & 97 & 87  &  83 & 82  &  81 &  79 & 72 & 67\\
grow15    &113 & 97 &101&94 & 88& 81  & 71& 66 \\
grow22    &112&109&103&89&82&72& 63&60 \\
israel    & 54  & 55   &67  & 82   & 91  & 118  &  153 &$\star\star$\\
kb2       & 38  & 42   &  36  & 40   &  45  & 61   & 65   & $\star\star$  \\
lotfi     & 46  & 46   &  48  & 53   & 58  & 66   & 88   & 106  \\
maros     & 58  & 61   &  63  &66 & 77  & 101  &164  & $\star\star$  \\
maros-r7  & 30  & 29   &  29  & 29   &  30  & 32   & 37   & 48 \\
modszk1   & $\star\star$ & 58   &  55  & 59   &58  & 75   & 80   & $\star\star$ \\
nesm      & 97  & 97   & 92  & 98   &  104 & 115&$\star\star$&$\star\star$\\
perold    &$\star\star$&$\star\star$&116  &146 &198&288&$\star\star$&$\star\star$\\
pilot     &100 &$\star\star$ &101& 108 &121 & $\star\star$&$\star\star$&$\star\star$\\
pilot4    & 125&144&144 & 173& 219  &$\star\star$&$\star\star$&$\star\star$\\
pilot87   & 116&118&121& 136& 178 &247&$\star\star$&$\star\star$\\
pilot\_ja  & $\star\star$&162&265&$\star\star$& $\star\star$ & $227^*$ &$\star\star$&$\star\star$\\
pilotnov  & 47  & 53   &  58  & 70   &  87  & 129 &$\star\star$&$\star\star$\\
pilot\_we  &$\star\star$&$\star\star$&$\star\star$&$\star\star$&$229^*$&229&$\star\star$&$\star\star$\\
qap8      & $39$  & $36$   &  $36$  & 32  &  34  & 37   & 37    & 38\\
qap12     & 59  & 53   &  55  & 54   & 53  & 53   & 55    &54\\
qap15     & $219^*$ &$168^*$ &  66 & 61   &58 &61   & $191^*$    & 68\\
recipe    & 38  & 38   &  37  & 38   &  38  & 33   & 31    & 48\\
sc105     & 26  & 27   &  27  & 28   &  26  & 29   & 29    & 32\\
sc205     & 35  & 37   &  31  & 34   &  35  & 40   & 47    & 63\\
sc50a     & 33  & 32   &  23  & 23   &  23  & 22   & 23    & 25\\
sc50b     & 23  & 23   & 22   & 21   &  21  & 21   & 20    & 21\\
scagr25   & 31  & 32   &  32 & 34   &  36  & 39   & 45    & 63\\
scagr7    & 37  & 29   &  31 & 31   &  32  & 42   & 39    & 52\\
scfxm1    & 44  & 44   &  45  & 48   &  53  & 61   & 83    & 118\\
scfxm2    & 46  & 47   &  51  & 54   &  62  & 74   & 95    & 148\\
scfxm3    & 45  & 45   &  50 & 54   &  64  & 81  & 102    & 146\\
scorpion  & 45  & 44   &  46  & 44   &  46  & $\star\star$ & $\star\star$&$\star\star$\\
scrs8     &$\star\star$&$\star\star$&$\star\star$&$\star\star$& 107& 134 & 190 &$\star\star$\\
scsd1     & 47  &45   & 44  & 41  & 41  & 37   & 40  & 41 \\
scsd6     & 53  & 47   & 46 &  43  &  41  & 41   &  32   & 34 \\
scsd8     & 42  & 41   &  39  &  38  &  37  & 37   &  29   & 30\\
sctap1    & 49  & 49   &  51  &  51  &  55  & 52   &  62   & 73 \\
sctap2    & $52$  & 51   &  48  &  48  &  49  & 46   &  54   & 67 \\
sctap3    & 56  & 49   &  47  &  46  &  48  & 50   &  59   & 73 \\
seba      & 47  & 48   &  47  &  43  &  45  & 56   &  62   &99\\
share1b   & 59  & 66   & 73   & 113  & 119  & 145  & 141   & $\star\star$\\
share2b   & 37  & 39   & 38   &  38  &  29  & 33   &  33   & 50\\
shell     & 51  & 48   &  47  &  45  &  46  & 48   &  54   & 67\\
ship04l   & 52  & 50   &  49  &  47  &  46  & 37   &  37   & 42\\
ship04s   & 52  & 50   &  48  &  46  &  36  & 45   &  35   & 39\\
ship08l   & 49  & 48   &  48  &  45  &  36  & 37   &  34   & 38\\
ship08s   & 52  & 50   &  49 &   48  &  46  & 39   &  37   & 50 \\
ship12l   & 51  & 49   &  48  &  46  &  44  & 45   &  49   & $\star\star$ \\
ship12s   & 49  & 48   &  46  &  46  &  40  & 40   &  49   & 53 \\
sierra    & 52  & 51   &  51  &  54  &  68  & 86   &  141  & $\star\star$\\
stair     & 46  & 42   &  34  &  45  &  41  & 46   &  59  & 99\\
standata  & 56  & 63   &  57  &  57  &  63  & 63   &  75   & 89 \\
standgub  & 56  & 54   &  53  &  57  &  55  & 58   &  73   & 88\\
standmps  & 71  & 66   &  64 &   65  &  68  & 69   &  81   &101\\
stocfor1  & 43  & 43   &  43  &  45  &  49  & 48   &  56   & 73\\
stocfor2  & 70  & 72   &  78  &  89  &  99  &116   & $162^*$   & $232^*$\\
stocfor3  & 49  & 46   &  47  &  46  &  47  & 48   & 51    & 56\\
truss     & 49  & 46   &  47  &  46  &  47  & 48   &  51   & 56 \\
tuff      & 50  & 47   &  45  &  44  &  53  & 52   &  57   & 75\\
vtp.base  & 46  & 46   &  50  &  58  &  73  & 94   &  90   & $\star\star$\\
wood1p    &$\star\star$ & 76 &   63  &  60  &$\star\star$ & $\star\star$ &50&41\\
woodw     & 79  & 71   &  64  &  61  &$\star\star$ & $\star\star$&$\star\star$&$\star\star$ \\ \hline
\end{longtable}
\centerline{\bf Table 1}

$\star\star$ : Number of iterations exceeds 300.

$*$: Best optimal value obtained with $Rgap\in(10^{-8},10^{-10})$.

At first we can observe that for every problem there is at least one $r$ value for which the problem is solved. Also, we can see that most problems are solved for $r$ between 0 and 0.5. The maximum number of problems solved is reached for $r = 0.2$, as shown in the graph below.

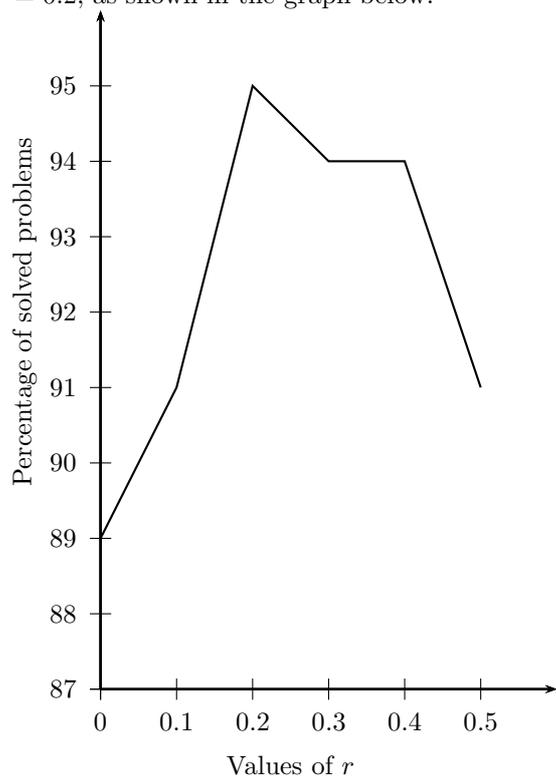
\centerline{
\begin{pspicture}(0,87)(8,96)
\psaxes[Ox=0,Oy=87,Dx=0.1,dx=1,Dy=1]{->}(0,87)(6,96)
\rput{90}(-1,92){Percentage of solved problems}
\rput(2.5,86){Values of $r$}
\psline(0,89)(1,91)(2,95)(3,94)(4,94)(5,91)
\end{pspicture}
}


\section{Concluding remarks}\label{remarks}

The results are conclusive and show that differentially barriers penalty functions offer effective alternative to conventional logarithmic barrier function in linear programming.

As we point out in the introduction, we chose an algorithm of affine scaling type for the simplicity of its implementation. But the "Predictor-corrector" method of Mehrotra \cite{meh} has proved highly efficient in the classical case ($r = 0$). Our immediate goal is to adapt this method to these new penalty functions.


\end{document}